\newcommand\AGACSE{0}  
\ifnum0=\AGACSE
\documentclass[12pt]{article}

\usepackage[english]{babel}

\usepackage[a4paper,top=2cm,bottom=2.5cm]{geometry}
\usepackage[colorlinks=true, allcolors=blue]{hyperref}
\newcommand{\inst}[1]{}
\newcommand{\orcidID}[1]{}
\newcommand{\titlerunning}[1]{}
\newcommand{\authorrunning}[1]{}
\newcommand{\institute}[1]{\date{#1}}
\newcommand{\email}[1]{}
\newcommand{\Keywords}[1]{\par\vspace{1ex}\noindent
        {\bf Keywords} #1}

\else 
\documentclass{jssc}
\usepackage{hyperref}
\usepackage{color}

\input jsscN.tex
\fi

\usepackage{amsmath}
\usepackage{amssymb}
\usepackage{amsthm}
\usepackage{multirow}
\usepackage{graphicx,subcaption}

\usepackage{anyfontsize}
\usepackage[T1]{fontenc} 
\usepackage{diagbox}
\newtheorem*{corollaryX}{Corollary}
\newcommand\thr[1]{Theorem~\ref{thm:#1}}
\newcommand\ler[1]{Lemma~\ref{lem:#1}}
\ifnum0=\AGACSE
\newtheorem{theorem}{Theorem}
\newtheorem{lemma}{Lemma}
\theoremstyle{definition}
\newtheorem*{definition}{Definition}
\newcommand\tab[3]{\begin{table}#1\caption{\label{tab:#2}#3}\centering}
\newcommand\tae{\end{table}}
\newcommand\fib[3]{\begin{figure}#1\begin{center}
   } 
\newcommand\fie[2]{\end{center}\vskip -4ex
    \caption{\label{fig:#1}#2}\end{figure}}
\newcommand\excessKeywords{, octonions, sedenions}
\else
\newcommand\tab[3]{\begin{table}\begin{minipage}{\textwidth}%
  \centering\caption{\label{tab:#2}#3}}
\newcommand\tae{\end{minipage}\end{table}}
\newcommand\fib[3]{\begin{figure}\begin{center}
    \includegraphics#2{figures/#3}
    }
\newcommand\fie[2]{\end{center}\vskip -2ex
    \caption{\label{fig:#1}#2}\end{figure}}
\newcommand\excessKeywords{}
\fi
\newcommand\Z{{\mathbb Z}} 
\newcommand\N{{\mathbb N}} 
\newcommand\R{{\mathbb R}} 
\newcommand\C{{\mathbb C}} 
\newcommand\HB{{\mathbb H}} 
\newcommand\OB{{\mathbb O}} 
\newcommand\SB{{\mathbb S}} 
\newcommand\PB[1]{{\mathbb P}_{#1}} 
\newcommand\ZB[1]{{\mathbb Z}_{#1}}
\newcommand\UB[1]{{\mathbb U_{#1}}}

\renewcommand\o[1]{{\rm o}_{#1}}
\newcommand\GT{{\rm G_2}}
\newcommand\Cl[1]{{\rm Cl}({#1})}
\newcommand\Cp[1]{{\rm Cl}^+({#1})}
\newcommand\Ci[1]{{\rm Cl}^*({#1})}
\newcommand\Aut[1]{{\rm Aut}({#1})}
\newcommand\AB[1]{{\rm A_{#1}}}

\newcommand\Sym[1]{{\rm S}_{#1}}
\newcommand\Spin[1]{{\rm Spin}({#1})}
\newcommand\Pin[1]{{\rm Pin}({#1})}
\newcommand{\Sharp}[1]{{\rm Sharp}({#1})}
\newcommand{\Harp}[1]{{\rm Harp}({#1})}
\newcommand{\Shar}[2]{{\rm Sharp}_{#1}({#2})}
\newcommand{\Har}[2]{{\rm Harp}_{#1}({#2})}
\newcommand{\PG}[2]{{\rm PG}(#1,#2)}
\newcommand\Rot[1]{{\rm R}_{#1}} 
\newcommand\degs{^\circ}

\newcommand\GG[2]{{\rm {#1}^{(#2)}}}
\newcommand\GA[2]{{\rm A}_{#1}^{(#2)}}
\renewcommand\o[1]{{\rm o_{#1}}} 
\newcommand\oo[1]{{\rm o_{#1}}} 
 
\newcommand\e[1]{{\rm e_{#1}}} 
\newcommand\ee[1]{{\rm e}_{#1}}

\newcommand\eqb[1]{\begin{equation}\label{eqn:#1}}
\newcommand\eqe{\end{equation}}
\newcommand\eqn[1]{{\rm (\ref{eqn:#1})}}
\newcommand\eab[1]{\[\begin{array}{#1}}
\newcommand\eae{\end{array}\]}
\newcommand\elb[1]{\begin{equation}\label{eqn:#1}\begin{aligned}}
\newcommand\ele{\end{aligned}\end{equation}}
\newcommand\tar[1]{Table~\ref{tab:#1}}
\newcommand\tas[2]{Tables~\ref{tab:#1}--\ref{tab:#2}}
\newcommand\fig[1]{Figure~\ref{fig:#1}}

\newcommand\ec{\text{,}}
\newcommand\es{\text{.}}
\newcommand{\sstrut}{\rule{0pt}{2.1ex}}
\newcommand{\tstrut}{\rule{0pt}{2.8ex}}
\newcommand{\hod}{^*}
\newcommand\nz{$\star$}
\newcommand\maps{:~}
\newcommand\mapa[2]{\langle {#1} :\:{#2}\rangle}
\newcommand{\smp}{{\scriptstyle +}}
\newcommand{\ssmp}{{\scriptscriptstyle +}}
\DeclareMathSymbol{\smu}{\mathbin}{AMSa}{"39}
\newcommand{\sm}{\texorpdfstring{\smu}{-}}
\def\abs{\mathop{\rm abs}\nolimits}

\begin{document}
\ifnum0=\AGACSE
\title{Automorphisms of Sedenions}
\titlerunning{Automorphisms of Sedenions}
\author{G.\:P.\:Wilmot\inst{}\orcidID{0009-0003-7451-3933}}
\authorrunning{G.\:P.\:Wilmot}
\institute{University of Adelaide, Adelaide, South Australia, 5005,
   Australia, \url{https://www.adelaide.edu.au}}
\maketitle\thispagestyle{empty}
\else
\title{Automorphisms of Sedenions$^*$}
{\uppercase{Wilmot} G.~P.}
{University of Adelaide, Adelaide, South Australia, 5005, Australia, 
  \url{https://www.adelaide.edu.au}} 
{$^*$This research was supported by an Australian Government Research Training Program Scholarship and an Australian Research Council Grant (FL240100217).}
\drd{DOI: 10.48550/arXiv.2512.07210}{Received: Dec 2025}{}
\dshm{2026}{XX}{Automorphisms of Sedenions}{\uppercase{Wilmot} G.~P.}
\fi

\newcommand\AbstractA{This paper extends the octonion calibration in Clifford algebra to two related sedenion calibrations. Associative calibrations map quaternion rings in Clifford algebra to those in Cayley-Dickson algebras, with octonions consisting of seven rings and sedenions having 35 rings. A new non-associative calibration is found that is related to the sedenion associative calibration but is invertible, providing a classification of the ideals of Clifford algebra and hence of $8\times8$ matrices.}
\newcommand\AbstractB{The non-associative calibrations generate an infinite series of new Clifford subalgebras that match the finite geometry $\PG{N}2$ series, discovered by Gino Fano, and connect the projections of simplices with those from hyper-cubes.
The visual representation of sedenions as the $\PG{3}2$ tetrahedron, representing the 15 Fano planes of sedenions, allows the sedenion power-associative loops to be visually distinguished from the octonion loops. The visual representation of sedenions may be important for particle physics and Clifford algebra provides analysis tools for the $\PG{N}2$ series that uncovers the hidden structure of Cayley-Dickson algebras within an associative algebra.}
\ifnum0=\AGACSE
\begin{abstract}\AbstractA
This leads to new subalgebras of Clifford algebras and a thorough analysis of the possible automorphisms of sedenions is applied using the Spin and Pin groups. 
While the calibrations and automorphisms of octonions have minimal divergence, it is found sedenions extend Lie exceptional algebra, $\GT$, to a Lie algebra of rank 49. This work clarifies and extends Schafer's and Brown's sedenion automorphisms.
Brown introduced automorphisms that change the sedenion representation, unlike $\GT$, but Clifford algebra analysis finds there are 92 such pure representations.
\AbstractB
\end{abstract}
\else
\Abstract{\AbstractA
While the calibrations and automorphisms of octonions have minimal divergence, it is found sedenions extend the exceptional Lie algebra $\GT$ to a Lie algebra of rank 49. This work clarifies and extends Schafer's and Browns's sedenion automorphisms.\AbstractB}
\fi

\section{Introduction}
The series of Cayley-Dickson algebras, $\AB{N}$, provides $n=2^N$ basis elements generated by $N$ generators at each level, $N\ge0$. The series starts with the reals, $\AB0\cong\R$, complex numbers, $\AB1\cong\C$, quaternions, $\AB2\cong\HB$, octonions, $\AB3\cong\OB$, followed by the sedenions, $\AB4\cong\SB$. The series has increasingly complicated norms until sedenions and all $\AB{N}$, for $N\ge4$, which have zero divisors, hence they are power-associative and cannot be normed algebras. Because of this, until recently, they have been avoided by physicists as not being appropriate for reality. Mathematicians have struggled with the analysis of these algebras concentrating on the alternate rule, which is only valid for octonions~\cite{Morano,Zhilina,Wilmot2}, and insufficient for embeddings of Lie algebras in the automorphisms of $\AB{N}$, for $N\ge4$.

The Lie algebras were formulated and classified by Killing and Cartan~\cite{Dray} and the smallest exceptional Lie algebra, $\GT$, was shown to be the automorphisms of octonions simultaneously by Cartan and Engel, $\Aut\OB\cong\GT$~\cite{Agricola}. Cartan found it as the solution space of a system of second order differential equations and Engel used contact transformations of Grassmannians, which fell out of favour but is the essence of the Clifford algebra approach. Schafer~\cite{Schafer} found that $\Aut{\SB}\overset{{?}}{\cong}\GT$ while, subsequently, Brown~\cite{Brown} found that $\Aut{\SB}\overset{{?}}{\cong}\GT\times\Sym3$. This discrepancy is analysed extensively in this paper with the main result being a new calibration that provides for the enumeration of the automorphisms of sedenions. This analysis extends the work of Schafer~\cite{Schafer} and Brown~\cite{Brown} for sedenions to find that the full set of automorphisms are more complicated than a direct product with $\GT$. This result disagrees with Brown's full result that $\Aut{\AB{N}}\cong\Aut{\AB{N-1}}\times\Sym3$ for $N\in(4,5,6)$ and the procedure used herein should be applicable for all higher order automorphisms of Cayley-Dickson algebras.

A graded notation is used for both Clifford and Cayley-Dickson algebras, but for different reasons. The basis of the $\Cl{n}$ algebras naturally represent the same geometric objects as Grassmann's exterior algebra and the grading follows that of differential geometry. The $\AB{N}$ algebras are not directly related to geometry but the doubling rule that defines the algebra~\cite{Harvey} introduces a new generator at each level with the binomial expansion of all generators providing the grading, as explained in the following.

\begin{definition}
The basis of the Clifford algebra~\cite{Harvey} of positive definite signature is composed of basis elements $\Cl{n} = \langle\e{\alpha}\maps\alpha\in\{\emptyset,1,2,\dots,n,12,13,\dots,12\ldots n\}\rangle$. Clifford algebra is the real module of the Pin group, $\Cl{n} = \langle\e{\alpha}+\e{\beta}+\cdots\maps\e\alpha,\e\beta,\cdots\in\Pin{n}\rangle$. Any basis element, $\e\alpha$, can be mapped to Grassmann's exterior algebra and are called forms. The even subalgebra, $\Cp{n}$, is thus the real module of the $\Spin{n}$ group. The basis of the Cayley-Dickson algebra is denoted as $\AB{N} = \langle\o{\alpha}\maps\alpha\in\{\emptyset,1,2,\dots,n,12,13,\dots,12\ldots n\}\rangle$ but has a different structure due to multiplication being defined by the doubling rule~\cite{Harvey}. Unity is used in these algebras as $\e{\emptyset}=\o{\emptyset}=1$ and $\ee{\alpha}^2=1$ and $\oo{\alpha}^2=-1$ for $\alpha\neq\emptyset$, using Porteous's notation~\cite{Porteous} for a subset of natural numbers. The grading of these three algebras exposes the XOR multiplication operation, $\alpha\veebar\beta$, but the resulting sign and ordering is different for the Cayley-Dickson algebra. The ordering for $\Cl{n}$ is the same as the geometry of Grassmann's algebra but $\AB{N}$ uses the binomial theorem doubling rule~\cite{Wilmot2}. 
\end{definition}
Some of the sub-rings of $\Cl{n}$, such as quaternions, form groups. Here we wish to talk about ideals so it is important to distinguish the two quaternion rings, quaternions or anti-quaternions, depending on whether their triple product is $-1$ or $+1$, respectively. This paper classifies the $\Cp{n}$ idempotents, which highlights the sub-ring structure. Since these idempotents are derived from elements that are invertible, this work extends the invertible Clifford group, $\Ci{n}$, beyond the $\Pin{n}$ group.
\begin{definition}
A {\it calibration}~\cite{Harvey}, $\Phi$, in differential geometry is a $p$-form on a subset, $U$, of $\R^n$ satisfying $d\Phi=0$ on $U$ and has maximum one for any point in $U$ on the Grassmannian of oriented $p$-dimensional subspaces of $\Lambda\R^n$. An associative calibration is a 3-form that corresponds to associative algebras where each term represents a quaternion group.
\end{definition}
The 3-form is the calibration derived from $\e{123}={\mathrm d}(x_1\e{23} -x_2\e{13} +x_3\e{12})/3$. In Clifford algebra $\e{123}$ is an element of $\Pin3$, which is represented by the face of the 2-simplex. The 2-forms of $\Spin3$ are isomorphic to quaternions in $\Cl3$ with the mapping to Cayley-Dickson quaternions, $\HB$,
\[ (\e{23},\e{13},\e{12})\rightarrow (\o1,\o2,\o{12})\ec \]
The parity of Clifford-algebra allows mappings to be defined to the quaternions and anti-quaternion rings, as
\elb{map3}
    \e{123}&\maps\{(\e1,\e2,\e3)\rightarrow (\o1,\o2,\o{12})\}\text{ and}\\
    \e{321}&\maps\{(\e1,\e2,\e3)\rightarrow (\o{12},\o2,\o1)\}\ec
\ele
where $\o1^2=\o2^2=\o{12}^2=-1$, $\o1\o2=\o{12}$, $\o1\o2\o{12}=-1$ and  $\o{12}\o2\o1=1$.

\begin{definition}
{\it Associative calibrations} define faces of the $(n-1)$-simplex that only touch each other at vertices and cover all edges, which is a representation of independent quaternions in $\Cl{n}$ that map to quaternions in Cayley-Dickson algebras. These define cross products and thus a multiplication table that corresponds to Cayley-Dickson algebras.
\end{definition}
\begin{definition}
Each associative calibration has many representations but only one corresponds to the faces mapping to only quaternion rings in the Cayley-Dickson algebra. Thus there is a concept of {\it pure calibrations} and only a single index is needed for each calibration representation.
\end{definition}
The first pure seven dimensional associative calibration~\cite{Harvey,Bryant,Wilmot1}, is
\eqb{cal7} \theta_1 = \e{123}+\e{145}-\e{167}+\e{246}+\e{257}+\e{347}-\e{356}\es\eqe
There are seven reflections of $\theta_1$, with eight negations, that are also calibrations that do not map to seven quaternions but include several anti-quaternions. There are 30 pure calibrations, $\theta_i$, $i\in\N_1^{30}$, with maps equivalent to the following for $\theta_1$,
\eqb{map7}
  \e{1234567}\maps\{(\e1,\e2,\e3,\e4,\e5,\e6,\e7)\longrightarrow
          (\o1,\o2,\o{12},\o3,\o{13},\o{23},\o{123})\}\es
\eqe
The two negative terms in \eqn{cal7} map to anti-quaternions, $\o1\o{23}=-\o{123}$ and $\o{12}\o{13}=-\o{23}$, and the parity, $-\e{167}$ and $-\e{356}$, change these to quaternion rings.
\begin{definition}
The {\it dual Clifford algebra} is composed of basis terms, $\e\alpha$, multiplied by the negated pseudoscalar and denoted $\e\alpha\hod$.
\end{definition}

In $\Cl7$, the pseudoscalar is $\e{1234567}$ and the dual of the calibration $\theta_1$ defines the co-associative calibration, 
\elb{dcal7} \theta_1\hod &= -\e{1234567}\theta_1 \\
       &= -\e{1247}+\e{1256}+\e{1346}+\e{1357}-\e{2345}+\e{2367}+\e{4567}\es
\ele
It was shown in~\cite{Wilmot1} that invertible forms can be created from the calibration $\theta_i$, $i\in\N_1^{30}$, and its dual in $\Cl7$, for $k$ the parity of $\theta_i$ or the number of minus signs, as
\elb{inv7}
  (3\e{1234567} -(-1)^{k}\theta_i)^2 &= -16\text{ and} \\
  (3 -(-1)^{k}\theta_i\hod)^2 &= 16\es
\ele
Hence $\theta_1(6\e{1234567}-\theta_1) = 7$.
\begin{definition}
The terms of a calibration with signs that are not related to a form that is invertible will be called a {\it quasi-calibration}. Calibrations and quasi-calibrations are together called {\it generalised calibrations}. A generalised calibration with all positive terms will be called a {\it primary calibration} or just primary.
\end{definition}
Thus there are 30 primary and pure calibrations in $\Cl7$ with $2^{7}=128$ generalised calibrations each, 16 of which are calibrations giving 480 in total.

The quasi-calibrations generate other algebras called quasi-octonions. Out of the $2^7\times30-480$ possible algebras,~\cite{Wilmot1} showed that only six are generated, denoted as $\PB{k}$ for $k\in(4,8,10,12,14,16)$. The first, $\PB4$, was found by Cawagas~\cite{Cawagas1} and~\cite{Wilmot2} proved that only three quasi-octonions are included in all Cayley-Dickson algebras. The embeddings of these quasi-octonion algebras are shown in \fig{ultron} where $\UB{m}\hod=\AB{m+3}$, $m>0$, and the notation, $\UB{m}$, $m\ge0$, called ultracomplex numbers, denotes all power-associative algebras including the quasi-algebras and power-associative Cayley-Dickson algebras, $\AB{N}$, $N>3$. The notation $\UB0$ denotes all six quasi-octonions but generally not octonions, since they are not power-associative.
\fib{[ht]}{}{}\includegraphics[width=10cm]{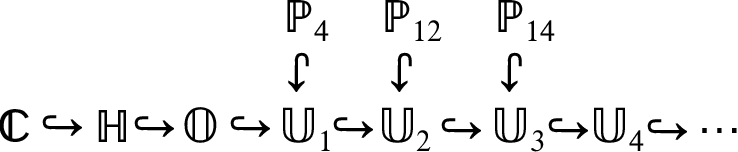}
\fie{ultron}{Cayley-Dickson Algebras Structure}

\begin{definition}
A {\it triad} refers to a triple of Cayley-Dickson basis elements that are distinct and pure, or non-scalar.
The {\it ordered associator} for a triad is defined as $[a,b,c]=(ab)c-a(bc)$, for $a<b<c$ in the ordering provided by the grading~\cite{Wilmot2}. If $[a,c,b]-[c,b,a]+[b,a,c]=0$ then the triad is called {\it associative}. Quaternions and anti-quaternions are associative.
\end{definition}

\begin{definition} 
The associator is extended in~\cite{Wilmot2} to an {\it unordered associator} that defines four non-associative types, shown below for the triad with $a<b<c$, where \nz\  denotes a non-zero value.
\begin{center} \begin{tabular}{|c|c|c|c|} \hline
\multicolumn{3}{|c|}{\bf Associativity} 
   &\multirow{2}{*}{\bf Non-Associativity} \\\cline{1-3}
 $[b,a,c]$ &$[a,b,c]$ &$[a,c,b]$ &\\\hline
 \nz  &0   &0   &{\bf Type~A}\\\hline
 0    &\nz &0   &{\bf Type~B}\\\hline
 0    &0   &\nz &{\bf Type~C}\\\hline
 \nz  &\nz &\nz &{\bf Type~X}\\\hline
\end{tabular} \end{center}
\end{definition}

The usual ordered associator only counts one of these columns, which is the middle one for ordered triads, so that only Types~B and~X are counted. But, of course, there are only seven quaternion rings in the quasi-octonions so the remainder of the 28 non-associative triads to make up the $\binom73=35$ faces of the 6-simplex are of Types~A and~C. It was found in~\cite{Wilmot1} that each quasi-octonion has 12 zero divisors and~\cite{Wilmot2} proved each quasi-octonion is embedded in factors of seven. Since the zero divisors come from the quasi-octonion rings, \fig{ultron} shows each Cayley-Dickson algebra after octonions has factors of 84 zero divisors.
\tab{[ht]}{pbk}{Algebra Identification}
\begin{tabular}{|c|l|l|l|l|l|}  \hline 
{\bf Algebra} &{\bf A} &{\bf B} &{\bf C} &{\bf X} &{\bf Generalised Calibration} \\\hline
$\OB$    &0 &0 &0&28 &$\e{123}+\e{145}+\e{167}-\e{246}+\e{257}+\e{347}+\e{356}$\\\hline
$\PB4$   &12&0 &12&4 &$\e{123}+\e{145}+\e{167}+\e{246}-\e{257}+\e{347}+\e{356}$\\\hline
$\PB8$   &10&4 &10&4 &$\e{124}+\e{135}+\e{167}+\e{237}+\e{256}+\e{346}+\e{457}$\\\hline
$\PB{10}$&9 &6 &9 &4 &$\e{124}+\e{135}+\e{167}+\e{237}+\e{256}-\e{346}+\e{457}$\\\hline
$\PB{12}$&8 &8 &8 &4 &$\e{123}+\e{145}+\e{167}+\e{246}+\e{257}-\e{347}+\e{356}$\\\hline
$\PB{14}$&7 &10&7 &4 &$\e{123}+\e{145}+\e{167}+\e{246}+\e{257}+\e{347}+\e{356}$\\\hline
$\PB{16}$&6 &12&6 &4 &$\e{125}+\e{134}+\e{167}+\e{236}+\e{247}+\e{357}+\e{456}$\\\hline
\end{tabular}\tae

The $\PB{k}$ have increasing cardinality of Type~B as shown in \tar{pbk}, which also provides the first calibration or quasi-calibration for each algebra. It was shown in~\cite{Wilmot2} that all Type~A triads and half of Type~B generate zero-divisor pairs so \tar{pbk} shows that all $\PB{k}$ have 12 zero divisor pairs thus enabling power-associativity. This table is an attempt to justify why only three of the six quasi-octonions appear within the Cayley-Dickson series. It can be seen that $\PB{k}$ for $k\in\{4,12,14\}$ occur in the first primary calibration while the others appear in subsequent primary forms. Unfortunately, although $\PB{12}$ and $\PB{14}$ appear in all 30 representations, it is found that $\PB4$ is missing from the second group of six representations whereas $\PB8$ and $\PB{16}$ are missing from half of the representations and $\PB{10}$ is missing from just two of the groups of six found in the 30 associative 3-form representations. This implies the first primary is special but this explanation is pure speculation.

This paper extends this work to sedenions distinguishing calibrations in $\Cl{15}$ from representations than generate quasi-sedenions. It starts by defining the Sharp algebras that extend the quaternions and the commuting 7-D algebra to what is conjectured to be an infinite series. This is initially a geometric proof that shows the power of the geometric Clifford algebra. The analogue to the 7-D categorisation of the ideals and the derivation of $\GT$ are provided in $\Cl{15}$. This explains the automorphisms of sedenions derived by Schafer and Brown and new terminology is introduced to distinguish isomorphisms of sedenions to other representations of themselves that are not automorphisms. The automorphisms correspond to a single domain, whereas herein we show three domains of isomorphisms with one corresponding to the $\PB4$ rings inside sedenions. This is more complicated than the isomorphisms of octonions and it is speculated that the $\UB2$ and $\UB3$ levels shown in \fig{ultron} will involve $\PB{12}$ and $\PB{14}$ domains providing isomorphism structures that are much more complicated than for sedenions.

\section{Sharp and Harp algebras}
The $\Sharp{N}$ and $\Harp{N}$ algebras are subalgebras of Clifford algebras, $\Cp{n}$ and $\Cl{n}$, respectively, where $n=2^N-1$, $N>1$. These are defined formally in the next section and introduced geometrically first.

The 2-simplex is a triangle in 2-D that can be projected into a line segment in 1-D, called just a line. It has the perspective geometry that each edge has two vertices and each vertex has two edges. Clifford algebra adds $k$-forms as labels to every geometric element and a parity as arrows to 2-forms and above. $\Sharp2$ is represented by a directed 2-simplex with consistent parity so that it can be projected to a 1-D directed line, as shown in \fig{2sim}. Taking any vertex label and multiplying by another vertex in $\Cl3$ gives the connecting edge, which is the dual of the third vertex. In other words the dual of one vertex multiplied by the dual of another gives the dual of the third and by the perspective geometry this applies to all edges so forms a ring. This is the $\Sharp2$ algebra, either quaternion or anti-quaternion, and the dual algebra uses elements of $\Pin3$, which contains the calibrations $\pm\e{123}$. The quaternions are not members of $\Spin3$ since $\e{12}+\e{13}+\e{23}$ is not a product of vectors but are included in $\Sharp2$ as a subalgebra of $\Cp3$. 
\fib{}{}{}\includegraphics[width=8cm]{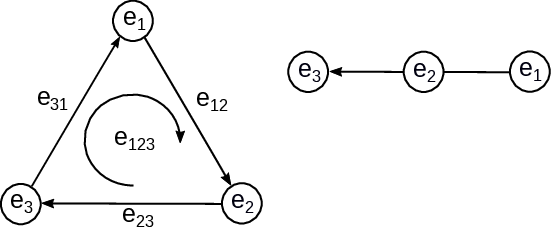}
\fie{2sim}{Directed 2-simplex and projection to 1-D}

The 6-simplex projects to the 2-D Fano plane shown in \fig{6sim} where each line represents a face with three edges that provides the $\binom72=21$ rotations in $\Cl7$. The 35 faces predicted by Pascal's triangle can also be counted. The perspective geometry is that each line intersects three vertices and each vertex joins three lines. The Clifford algebra interpretation is that multiplying the 3-form face provided by any line with the 3-form for another line contracts the common vertex to give the dual of the third intersecting line for the common vertex. Hence the dual of a line multiplied by the dual of another line gives the dual of the shared third line for the common vertex. By the projective geometry this is the 4-form $\Shar{i}3$ algebra provided by the terms of the co-associative calibration \eqn{dcal7} and including the dual terms gives the $\Har{i}3$ algebra, where $i\in\N_1^{30}$ denotes the calibration representation. Both \eqn{cal7} and \fig{6sim} specify $i=1$ in this case.
\fib{[ht]}{}{}\includegraphics[width=4.6cm]{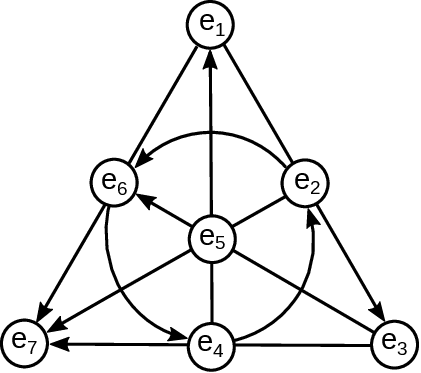}
\fie{6sim}{Directed Fano plane diagram}
\fib{[ht]}{}{}\includegraphics[width=9.7cm]{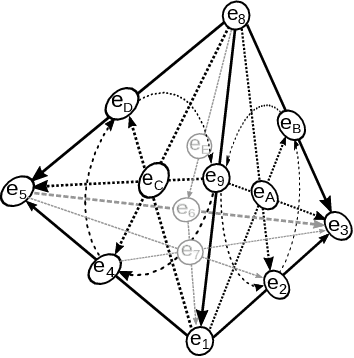}
\fie{14sim}{Fano volume diagram}

The 14-simplex projects to the 3-D tetrahedron shown in \fig{14sim}, where hexadecimal notation is used for the 15 dimensions of $\Cl{15}$. The 14-simplex contains 15 Fano planes whereby the projective geometry is that each plane shares a line with two other Fano planes. Also seven planes are incident on each vertex and each plane contains seven vertices. This is the perspective geometry Fano discovered in 1892 in his paper uncovering an infinite series of harmonic spaces related to quadrangles~\cite{Fano} and he gave the same formula for the number of quaternion lines and the number of Fano planes at each level, as provided in \thr{2}. Hence this is the same infinite series conjectured for simplex projections represented in Clifford algebra as Sharp algebras. The difference is that Fano deals with projections though a vertex of a hyper-cube whereas here we have projections from infinity of simplices. Since he ``imagined'' this series and preferred geometric proofs, the series of such figures should be known as the Fano hyper-volumes, starting with the Fano volume in 3-D, derived from the $\Sharp{4}$ algebra. 

Each vertex in \fig{14sim} intersects seven Fano planes and each plane joins seven vertices. There are 15 Fano planes and each line of $\Sharp{4}$ contains three edges, which gives the $15\times7=105=\binom{15}2$ edges of the 14-simplex. Also each quaternion line is shared between three planes and since the Fano planes represent a 7-form in $\Cl{15}$ then the example product of the two faces at the front of \fig{14sim} contracts the shared line, $\e{189}$, to result in the dual of the Fano plane from this line to the middle of the rear face, $\e{16789EF}$. By the perspective geometry, products of the duals of any two Fano planes results in the dual of the plane that shares the common quaternion line. The projective geometry is proved with the algebra because the product of any two 8-forms of $\Sharp{4}$ generates another member of the algebra. Like $\Sharp3$ for octonions, $\Sharp{4}$ could be called commuting, associative sedenions, since is has the same 3-cycle rings mapped from the quaternions.

Cawagas~\cite{Cawagas1,Cawagas2} found from loop analysis that eight octonions and seven $\PB4$ rings are embedded in the sedenions. This means the count of ordered non-associative triads matches the count for each subring, $252 = 8\times28 +7\times4$. These rings, provided here in reduced triad form that generate the 7-cycles, were found in~\cite{Wilmot2} and are replicated here as $\phi_O$ and $\phi_P$, representing octonion and $\PB4$ representations, respectively,
\eqb{phi}\begin{split}
  \phi_O = &((\o1,\o2,\o3),(\o1,\o2,\o4),(\o1,\o3,\o4),(\o1,\o{23},\o4),\\
  &(\o2,\o3,\o4),(\o2,\o{13},\o4),(\o{12},\o3,\o4),(\o{12},\o{13},\o4))\text{ and}\\
  \phi_P = &((\o1,\o2,\o{34}),(\o1,\o3,\o{24}),(\o1,\o{23},\o{24}),(\o2,\o3,\o{14}),\\
  &(\o2,\o{13},\o{14}),(\o{12},\o3,\o{14}),(\o{12},\o{13},\o{14})).
\end{split}\eqe

The inverse map applied to the triads of $\phi_O$ and $\phi_P$ can be used to define 7-forms by first expanding the triads in \eqn{phi} to the seven terms of each loop. For example, $(\o1,\o2,\o3)\rightarrow(\o1,\o2,\o{12},\o3,\o{13},\o{23},\o{123})$, which is inverse mapped by \eqn{map7} to $\e{1234567}$. This is the Fano plane on the bottom of \fig{14sim}, denoted $\Phi_A$, and represents the 7-dimensional octonion in the top left-hand quarter of sedenion products in \tar{sedmul}. The equivalent process applied to the remaining terms provides the following.
\begin{definition}
The non-associative rings of sedenions are reversed mapped to 7-forms in $\Cl{15}$, accumulated as a 7-form, $\Phi$, with three domains,
\eqb{cal15} \Phi=\Phi_A +\Phi_O +\Phi_P\es \eqe
The first $\phi_O$ term is recognised as a representation of the 7-D octonions, and the other domains are the remaining $\OB$ and $\PB4$ parts, 
\eqb{Phi}\begin{split}
\Phi_A &= \e{1234567}\ec \\
\Phi_O &= \e{12389AB} +\e{14589CD} +\e{3568BDE}\\
     &\quad +\e{2468ACE} +\e{2578ADF} +\e{3478BCF} +\e{16789EF}\text{ and}  \\
\Phi_P &= \e{123CDEF} +\e{145ABEF} +\e{3569ACF} \\
    &\quad +\e{2469BDF} +\e{2579BCE} +\e{3479ADE} +\e{167ABCD}\es\\
\end{split}\eqe
\end{definition}

Note that $\Phi_P = \e{89ABCDEF}\Phi_O$, which removes the apex vertex in \fig{14sim} from each $\Phi_O$ term and swaps indices from $\e{9ABCDEF}$, which is called hereinafter the false Fano plane. The term $\Phi_A$ is recognised as the bottom face of \fig{14sim} and the first three terms of $\Phi_O$ are the other three faces of the tetrahedron. The last three terms of $\Phi_O$ can be recognised as the three vertical Fano planes from median lines of $\Phi_A$ to the top vertex. And finally, there is the octonions cone, $\e{2468ACE}$, consisting of the bottom ring, the midpoint of the false Fano plane and the top vertex, $\e8$. This is shown in \fig{fanoc} and this projects to a 2-D Fano plane. 

\ifnum0=\AGACSE
\begin{figure}[!ht]\begin{center}
  \begin{subfigure}[b]{0.4\textwidth}
  \centering\includegraphics[width=4.7cm]{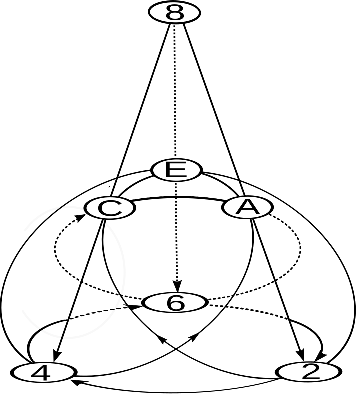}
  \caption{$\OB$ Fano Cone\label{fig:fanoc}}
  \end{subfigure}
  \begin{subfigure}[b]{0.4\textwidth}
  \centering\includegraphics[width=4.8cm]{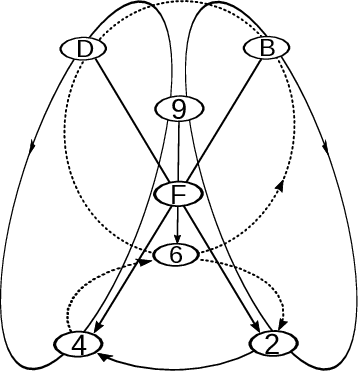}
  \caption{$\PB4$ Fano Folly \label{fig:fanof}}
  \end{subfigure}
  \vskip -1ex \caption{}
\end{center}\end{figure}
\else
\begin{figure}\begin{center}
  \subfigure[$\OB$ Fano Cone]{
  \centering\includegraphics[width=4.7cm]{figures/Article4-fig4.eps}
  \label{fig:fanoc}}
  \subfigure[$\PB4$ Fano Folly]{
  \centering\includegraphics[width=4.8cm]{figures/Article4-fig5.eps}
  \label{fig:fanof}}
  \vskip -1ex \caption{}
\end{center}\end{figure}
\fi
This leaves the $\PB4$ Fano planes from $\Phi_P$. The first three terms of $\Phi_O$ when removing $\e8$ and swapping terms of the false Fano plane rotate the outside faces down to the inner faces that share bottom edges and include the middle vertex, $\e{F}$, not shown in \fig{14sim}. The last three terms of $\Phi_O$ when swapping the upper indices swap the inner octonion planes with the middle terms to the edge face and middle terms of the other faces to give the cone to the adjacent vertex. For example, the inner plane to the right-hand face in \fig{14sim}, $\e{2578ADF}$, changes to $\e{2579BCE}$, the quasi-octonion cone from the right-hand face to the vertex $\e5$.

Finally, the remaining quasi-octonion Fano plane. Here we start with the octonion cone, $\e{2468ACE}$, and remove $\e8$ and swap the false Fano plane terms $\e{ACE}$ with $\e{9DBF}$. This is a thorn shaped diagram shown in \fig{fanof} that joins the four outside face rings and can be considered to hold the tetrahedron together. It is called the Fano folly because the directed simplices are Tits' building models of the orthogonal group and thus this is a Tits' apartment, which also projects to a Fano plane. The Fano volume represents the projective space $\PG32$~\cite{Polster}, where the folly is called the body centroid~\cite{Guy}. The mapping to sedenions from $\PG32$ is also known~\cite{Saniga} but the Clifford algebra labelling and the recognition of the $\PB4$ Fano planes in \fig{14sim} will show that Clifford algebra has hidden information about the non-associative rings of Cayley-Dickson algebras.

The graphical depiction of the Fano planes, cones and the folly in \fig{14sim} demonstrates the three fold embedding of quaternions in the octonion rings of sedenions.
Gresnigt~\cite{Gresnigt} uses such sharing of quaternions to represent coloured fermions but care must be taken to choose a quaternion ring not represented in $\phi_P$, as the zero divisors in one representation of a quark would make the family unsymmetrical. From \fig{14sim}, such quaternions that are triply shared by octonions must include $\e8$, which excludes the 7-dimensional Fano plane on the bottom of the Fano volume. Either an upper edge or the 7-D ring, $\e{246}$, is needed. One could speculate that the third Cayley-Dickson power-associative algebra, $\UB3$, which contains all three pseudo-octonion algebras, $\PB4$, $\PB{12}$ and $\PB{14}$, and look for quaternions spanning three power-associative subalgebras as a representation of three generations of coloured fermions.

\section{Non-Associative Calibrations}
\begin{definition}
Calibrations in $\Cl{n}$ are a sum of $k$-forms representing simplices whereby each vertex intersects $k$ lots of $(k-1)$-simplices and each of these simplices joins $k$ vertices. These are non-associative for $n>7$.
\end{definition}
We now show that $\Phi=\Phi_A +\Phi_O +\Phi_P$ is an invertible calibration and corresponds to the associative calibration that maps to sedenions. Starting with the 15-D associative calibration, rotated from the one found in~\cite{Wilmot1} as a representation of sedenions, $\SB$,
\elb{cal35}
 \Theta = \frac13&(\e{123} +\e{145} -\e{167} +\e{189} -\e{1AB} -\e{1CD} +\e{1EF}
    +\e{246} +\e{257}\\[-1ex]
    & +\e{28A} +\e{29B} -\e{2CE} -\e{2DF} +\e{347} -\e{356} +\e{38B} -\e{39A} -\e{3CF} \\
    & +\e{3DE} +\e{48C} +\e{49D} +\e{4AE} +\e{4BF} +\e{58D} -\e{59C} +\e{5AF} -\e{5BE}\\
    & +\e{68E} -\e{69F}  -\e{6AC} +\e{6BD} +\e{78F} +\e{79E} -\e{7AD} -\e{7BC})\ec
\ele
then a map can be defined from $\Cl{15}$ to $\SB$ that corresponds to this representation as
\eqb{map15}\begin{split}
    \e{123456789ABCDEF}\maps(\e1,\e2,\dots,\e{F}) \longrightarrow\: 
    &(\o1,\o2,\o{12},\o3,\o{13},\o{23},\o{123},\o4,\\
      &\ \o{14},\o{24},\o{124},\o{34},\o{134},\o{234},\o{1234})\es
\end{split}\eqe
This is pure since all 3-forms in $\Theta$ map to quaternion rings. For example, $\e{1AB}\rightarrow\o1\o{24}\ne\o{124}$ so that $\e{1AB}$ is negated in $\Theta$.

The 35 quaternions of $\Theta$ represent independent faces of the 14-simplex that cover all edges and touch only in vertices. Hence they define a cross product and specify the sedenion multiplication table shown in \tar{sedmul}. The seven terms in $\theta_1$ from $\Cl7$ generate the 7-dimensional octonions in the top left hand quarter of the table. The table uses the graded notation and is a different representation to that used by Cawagas~\cite{Cawagas2} who used the representation whereby the first octonion was the seventh sign variation of $\theta_1$.  As this is not a pure quaternion representation, $\Theta$ is preferred.

\tab{[!ht]}{sedmul}{Sedenion multiplication table in graded form}
{\renewcommand{\arraystretch}{1.4}\setlength\tabcolsep{2pt}\tiny
\begin{tabular}{|c||c|c|c|c|c|c|c|c|c|c|c|c|c|c|c|}\hline
       &$\o{1}$  &$\o{2}$  &$\o{12}$   &$\o{3}$  &$\o{13}$  &$\o{23}$  &$\o{123}$  &$\o{4}$  &$\o{14}$ 
          &$\o{24}$ &$\o{124}$  &$\o{34}$  &$\o{134}$ &$\o{234}$ &$\o{1234}$\\\hline
 $\o{1}$   &$\sm1$  &$\o{12}$  &$\sm\o{2}$  &$\o{13}$   &$\sm\o{3}$ &$\sm\o{123}$  &$\o{23}$ &$\o{14}$  
   &$\sm\o{4}$ &$\sm\o{124}$  &$\o{24}$  &$\sm\o{134}$  &$\o{34}$ &$\o{1234}$ &$\sm\o{234}$\\
 $\o{2}$   &$\sm\o{12}$   &$\sm1$  &$\o{1}$  &$\o{23}$  &$\o{123}$  &$\sm\o{3}$  &$\sm\o{13}$  &$\o{24}$
   &$\o{124}$  &$\sm\o{4}$  &$\sm\o{14}$  &$\sm\o{234}$ &$\sm\o{1234}$  &$\o{34}$ &$\o{134}$\\
 $\o{12}$   &$\o{2}$  &$\sm\o{1}$   &$\sm1$  &$\o{123}$  &$\sm\o{23}$  &$\o{13}$   &$\sm\o{3}$ &$\o{124}$ 
   &$\sm\o{24}$  &$\o{14}$  &$\sm\o{4}$ &$\sm\o{1234}$  &$\o{234}$ &$\sm\o{134}$  &$\o{34}$\\
 $\o{3}$   &$\sm\o{13}$  &$\sm\o{23}$ &$\sm\o{123}$   &$\sm1$   &$\o{1}$  &$\o{2}$  &$\o{12}$  &$\o{34}$
   &$\o{134}$ &$\o{234}$ &$\o{1234}$   &$\sm\o{4}$  &$\sm\o{14}$  &$\sm\o{24}$ &$\sm\o{124}$\\
 $\o{13}$   &$\o{3}$ &$\sm\o{123}$  &$\o{23}$   &$\sm\o{1}$   &$\sm1$  &$\sm\o{12}$   &$\o{2}$ &$\o{134}$
   &$\sm\o{34}$ &$\o{1234}$ &$\sm\o{234}$  &$\o{14}$   &$\sm\o{4}$ &$\o{124}$  &$\sm\o{24}$\\
 $\o{23}$  &$\o{123}$  &$\o{3}$  &$\sm\o{13}$   &$\sm\o{2}$  &$\o{12}$   &$\sm1$   &$\sm\o{1}$ &$\o{234}$
   &$\sm\o{1234}$  &$\sm\o{34}$ &$\o{134}$  &$\o{24}$  &$\sm\o{124}$  &$\sm\o{4}$  &$\o{14}$\\
 $\o{123}$   &$\sm\o{23}$  &$\o{13}$  &$\o{3}$  &$\sm\o{12}$   &$\sm\o{2}$  &$\o{1}$   &$\sm1$ &$\o{1234}$
   &$\o{234}$ &$\sm\o{134}$  &$\sm\o{34}$  &$\o{124}$  &$\o{24}$  &$\sm\o{14}$  &$\sm\o{4}$\\
 $\o{4}$   &$\sm\o{14}$  &$\sm\o{24}$ &$\sm\o{124}$  &$\sm\o{34}$  &$\sm\o{134}$ &$\sm\o{234}$
   &$\sm\o{1234}$   &$\sm1$   &$\o{1}$  &$\o{2}$  &$\o{12}$   &$\o{3}$  &$\o{13}$  &$\o{23}$ &$\o{123}$\\
 $\o{14}$   &$\o{4}$ &$\sm\o{124}$  &$\o{24}$  &$\sm\o{134}$  &$\o{34}$ &$\o{1234}$  &$\sm\o{234}$
   &$\sm\o{1}$   &$\sm1$  &$\sm\o{12}$  &$\o{2}$  &$\sm\o{13}$   &$\o{3}$ &$\o{123}$  &$\sm\o{23}$\\
 $\o{24}$  &$\o{124}$  &$\o{4}$  &$\sm\o{14}$  &$\sm\o{234}$ &$\sm\o{1234}$  &$\o{34}$  &$\o{134}$
   &$\sm\o{2}$  &$\o{12}$   &$\sm1$  &$\sm\o{1}$  &$\sm\o{23}$  &$\sm\o{123}$  &$\o{3}$  &$\o{13}$\\
 $\o{124}$   &$\sm\o{24}$  &$\o{14}$  &$\o{4}$ &$\sm\o{1234}$  &$\o{234}$ &$\sm\o{134}$  &$\o{34}$
   &$\sm\o{12}$   &$\sm\o{2}$  &$\o{1}$   &$\sm1$  &$\sm\o{123}$  &$\o{23}$  &$\sm\o{13}$  &$\o{3}$\\
 $\o{34}$  &$\o{134}$ &$\o{234}$ &$\o{1234}$   &$\o{4}$  &$\sm\o{14}$  &$\sm\o{24}$  &$\sm\o{124}$
   &$\sm\o{3}$  &$\o{13}$  &$\o{23}$ &$\o{123}$   &$\sm1$   &$\sm\o{1}$  &$\sm\o{2}$  &$\sm\o{12}$\\
 $\o{134}$   &$\sm\o{34}$ &$\o{1234}$ &$\sm\o{234}$  &$\o{14}$   &$\o{4}$ &$\o{124}$  &$\sm\o{24}$
   &$\sm\o{13}$   &$\sm\o{3}$ &$\o{123}$  &$\sm\o{23}$   &$\o{1}$   &$\sm1$  &$\o{12}$  &$\sm\o{2}$\\
 $\o{234}$  &$\sm\o{1234}$  &$\sm\o{34}$ &$\o{134}$  &$\o{24}$  &$\sm\o{124}$  &$\o{4}$  &$\o{14}$
   &$\sm\o{23}$  &$\sm\o{123}$  &$\sm\o{3}$  &$\o{13}$   &$\o{2}$  &$\sm\o{12}$   &$\sm1$  &$\o{1}$\\
 $\o{1234}$  &$\o{234}$ &$\sm\o{134}$  &$\sm\o{34}$  &$\o{124}$  &$\o{24}$  &$\sm\o{14}$   &$\o{4}$
   &$\sm\o{123}$  &$\o{23}$  &$\sm\o{13}$  &$\sm\o{3}$  &$\o{12}$   &$\o{2}$  &$\sm\o{1}$ &$\sm1$\\\hline
\end{tabular}}\vskip 1ex\tae
Any three indices of $\Phi_A$ or in terms of $\Phi_O$ mapped as basis elements by \eqn{map15} that are not quaternion generate octonions and any three mapped $\Phi_P$ term indices, not quaternion or X non-associative, generate the corresponding $\PB4$ representation. The breakdown of this structure is provided in \tar{fvo} which is related to the Fano planes in \fig{14sim}, as provided in the last column. There is a relationship between the two calibrations of sedenions, $\Theta$ and $\Phi$. \tar{fvo} shows the breakdown of $\Theta$ into the components that relate to $\Phi$ so that the sum of the components in the table is $9\Theta$.

\tab{[!ht]}{fvo}{Octonion and $\PB4$ representations}
\begin{tabular}{|c|l|l|c|} \hline
$i$ &{\bf $\Phi_i$\sstrut} &{\hfil\bf $\Theta_i$\hfil} &{\bf Class}\\\hline
$1$ &$\e{1234567}$\sstrut &$\e{123}+\e{145}-\e{167}+\e{246}+\e{257}+\e{347}-\e{356}$ &$\OB$ Face\\
$2$ &$\e{12389AB}$ &$\e{123}+\e{189}-\e{1AB}+\e{28A}+\e{29B}+\e{38B}-\e{39A}$ &$\OB$ Face\\
$3$ &$\e{14589CD}$ &$\e{145}+\e{189}-\e{1CD}+\e{48C}+\e{49D}+\e{58D}-\e{59C}$ &$\OB$ Face\\
$4$ &$\e{16789EF}$ &$-\e{167}+\e{189}+\e{1EF}+\e{68E}-\e{69F}+\e{78F}+\e{79E}$ &$\OB$ Plane\\
$5$ &$\e{2468ACE}$ &$\e{246}+\e{28A}-\e{2CE}+\e{48C}+\e{4AE}+\e{68E}-\e{6AC}$ &$\OB$ Cone\\
$6$ &$\e{2578ADF}$ &$\e{257}+\e{28A}-\e{2DF}+\e{58D}+\e{5AF}+\e{78F}-\e{7AD}$ &$\OB$ Plane\\
$7$ &$\e{3478BCF}$ &$\e{347}+\e{38B}-\e{3CF}+\e{48C}+\e{4BF}+\e{78F}-\e{7BC}$ &$\OB$ Plane\\
$8$ &$\e{3568BDE}$ &$-\e{356}+\e{38B}+\e{3DE}+\e{58D}-\e{5BE}+\e{68E}+\e{6BD}$ &$\OB$ Face\\
$9$ &$\e{123CDEF}$ &$\e{123}-\e{1CD}+\e{1EF}-\e{2CE}-\e{2DF}-\e{3CF}+\e{3DE}$ &$\PB4$ Plane\\
$10$ &$\e{145ABEF}$ &$\e{145}-\e{1AB}+\e{1EF}+\e{4AE}+\e{4BF}+\e{5AF}-\e{5BE}$ &$\PB4$ Plane\\
$11$ &$\e{167ABCD}$ &$-\e{167}-\e{1AB}-\e{1CD}-\e{6AC}+\e{6BD}-\e{7AD}-\e{7BC}$ &$\PB4$ Cone\\
$12$ &$\e{2469BDF}$ &$\e{246}+\e{29B}-\e{2DF}+\e{49D}+\e{4BF}-\e{69F}+\e{6BD}$ &$\PB4$ Folly\\
$13$ &$\e{2579BCE}$ &$\e{257}+\e{29B}-\e{2CE}-\e{59C}-\e{5BE}+\e{79E}-\e{7BC}$ &$\PB4$ Cone\\
$14$ &$\e{3479ADE}$ &$\e{347}-\e{39A}+\e{3DE}+\e{49D}+\e{4AE}+\e{79E}-\e{7AD}$ &$\PB4$ Cone\\
$15$ &$\e{3569ACF}$ &$-\e{356}-\e{39A}-\e{3CF}-\e{59C}+\e{5AF}-\e{69F}-\e{6AC}$ &$\PB4$ Plane\\\hline
\end{tabular}\tae
\begin{lemma}\label{lem:1}
Identifying the 15 terms of $\Phi$ as $\Phi_i$ with $\Phi_1=\Phi_A$, $i\in\N_2^8$ corresponding to the terms of $\Phi_O$ and $i\in\N_9^{15}$ corresponding to those of $\Phi_P$ and subsets of terms of $3\Theta$ labelled as $\Theta_i$, provided in \tar{fvo}, then the correspondence between the 7-form and 3-form $\Cl{15}$ calibrations is
\[\Theta_i^3 = \begin{cases}
    -43\Theta_i -42\Phi_i\ec         &\text{for $i\in\N_1^8$, or} \\
    -19\Theta_i +6\Phi_i+24\theta\ec &\text{for $i\in\N_9^{15}$.}
    \end{cases}\]
\end{lemma}
\begin{proof}
Each of the quaternions in $\Theta$ are embedded three times in $\Phi$ so that $\sum_{i=1}^{15}\Theta_i=9\Theta$. The 3-form octonions, $\Theta_i$ are invertible from \eqn{inv7} in the form $(3\Phi_i-\Theta_i)^2=-16$. The 3-form $\PB4$ representations are not invertible but have a remainder~\cite{Wilmot1} represented by $\theta$ above, which is one of the terms from $\Theta_i$.
Building the multiplication table for each 3-form and finding 28 X non-associative triads for $i\in\N_1^8$ or 4 X non-associative and 12 Type A, 12 type C and 4 Type X triads for $i\in\N_9^{15}$ identifies $\OB$ or $\PB4$ classes, respectively~\cite{Wilmot2}, as identified in \tar{pbk}. These correspond to the division in \tar{fvo}. 
\end{proof}

\begin{theorem}\label{thm:1}
The 7-form $\Phi=\Phi_A +\Phi_O +\Phi_P$, is a non-associative calibration with invertible form $\psi$, $\psi^2=-1$, where
  \eqb{inv15}\psi = \frac18(7\e{123456789ABCDEF} -\Phi)\es\eqe
The calibration $\Phi$ generates sedenions under the mapping \eqn{map15}. The index for the 15-D calibration representations are dropped due to the great number of calibrations.
\end{theorem}
\begin{proof}
The graphical presentation of the Fano volume shows $\Phi$ is a calibration and \tar{fvo} shows each quaterion appears three times and there are 15 Fano planes. Now, the $\Phi_O$ consists of the terms of $\theta_1$ multiplied by  the pseudoscalar without the 3-form term. This swaps $\Theta_i$ terms to the false Fano plane. For example, $\e{12389AB}$ becomes
\[\Theta_2'=\e{9AB} +\e{9CD} -\e{9EF} +\e{ACE} +\e{ADF} +\e{BCF} -\e{BDE}\es\]
This is a calibration and all false Fano plane parts of the $\Phi_O$ terms act as calibrations. Multiplying by terms of $\Phi_O$ acts on the terms of $\Theta_2$ and $\Theta_2'$ separately with imaginaries $\e{1234567}$ and $\e{9ABCDEF}$ so that multiplying the results by $\e{123456789ABCDEF}$ returns the terms to $\Phi_O$. The lower indices will have the same results as \eqn{cal7} with the false Fano planes terms having similar action. Hence, like \eqn{cal7}, $\Phi_O$ will have the invertible form
\eqb{obc}  (3\e{123456789ABCDEF} -\Phi_O)^2 = -16\es\eqe
This invertible form, $\Phi_O$, looking at each term in \eqn{Phi}, relates to $\Phi_P$ as
\eqb{obc1}
  \e{89ABCDEF}\Phi_O = \Phi_P\ec
\eqe
where $\e{1234567}\e{89ABCDEF}=\e{123456789ABCDEF}$ and $\e{1234567}$ commutes with $\Phi_O$ and $\Phi_P$, which means $\e{89ABCDEF}$ does too.
We start by expanding $\Phi_O$ in \eqn{obc} as
\begin{align*}
  \Phi_O^2 &= 9 +6\e{123456789ABCDEF}\Phi_O -16\\
        &= -7 +6\e{123456789ABCDEF}\Phi_O\es
\end{align*}
and substituting this and \eqn{obc1} into $(8\psi)^2$, gives
\begin{align*}
    64\psi^2 &= (7\e{123456789ABCDEF} -\e{1234567} -(1+\e{89ABCDEF})\Phi_O)^2\\
     &= -49 -1 +2(1+\e{89ABCDEF})\Phi_O^2 -14\e{123456789ABCDEF}(\e{1234567} \\
        &\quad+(1+\e{89ABCDEF})\Phi_O) +2\e{1234567}(1+\e{89ABCDEF})\Phi_O \\
     =& -50 -2(1+\e{89ABCDEF})(7-6\e{123456789ABCDEF}\Phi_O) \\
     &\quad+14\e{89ABCDEF} -12(\e{123456789ABCDEF}+\e{1234567})\Phi_O\\
     =& -50 -14 -14\e{89ABCDEF} +12(\e{123456789ABCDEF}+\e{1234567})\Phi_O \\
     &\quad+14\e{89ABCDEF} -12(\e{123456789ABCDEF}+\e{1234567})\Phi_O\\
     &= -64\ec
\end{align*}
so that $\psi^2=-1$. Finally, the map \eqn{map15} takes the first, second and fourth indices of each term of $\Phi$ to the appropriate triad in \eqn{phi} that were shown in~\cite{Wilmot2} to have X non-associativity for $\Phi_O$ terms and A non-associativity for $\Phi_P$ terms as provided in \ler{1}. Note that this means $\Phi(14\e{123456789ABCDEF}-\Phi)=15$.
\end{proof}
\begin{corollaryX}
  Any $90\degs$ rotation between $\OB$ Fano planes and $\PB4$ Fano planes can transform a sedenion calibration to a quasi-calibration. These are represented by terms
  $${\rm R}_Q = \mapa{\e{ij}}{i\in\N_1^7, j\in\N_9^F}\es$$
\end{corollaryX}
\begin{proof}
  Reference~\cite{Wilmot1} showed that octonion calibrations were invertible whereas quasi-octonion calibrations had a remainder. Since the only subalgebras are quaternion rings then the octonion loops shown in \tar{fvo} will transform to other calibrations under any of the 105 basis $90\degs$ rotations of $\Spin{15}$. This is not the case for the $\PB4$ subalgebras of sedenions that may change remainder or the calibration change to another representation whereby the remainder represents another category. Since the remainder is one of the terms of the calibration, $\Theta_i$, and their indices are derived from $\Phi_i$, $i\in\N_9^{15}$, then the effect of rotations given by ${\rm R}_Q$ on $\Phi_P$ only need to be considered. 

  The $\PB4$ 7-forms in \tar{fvo} each have three 7-D indices that occur thrice in each $\Phi_i$ and four false Fano-plane indices that each occur four times. Basis rotations that don't mix these planes, even those including $\e8$, will not change the representation of the calibration. But those from ${\rm R}_Q$ will mix upper and lower indices in an inconsistent way, changing the representation and sometimes having a different effect on the remainder. This is found to sometimes change the structure of the $\PB4$ plane to another quasi-calibration, including $\PB8$ and $\PB{10}$. For example, the Fano Folly, $\Theta_{12}$, rotated by $90\degs$ in the plane $\e{19}$ becomes $\e{12B} +\e{14D} -\e{16F} +\e{246} -\e{2DF} +\e{4BF} +\e{6BD}$, which is a quasi-octonion that maps to $\PB{10}$. Such cases are found in~\cite{Wilmot2} to not be representations of sedenions. This is observed to occur in at least two $\PB4$ planes of \tar{fvo} for each of ${\rm R}_Q$. This is a general principle since of the eight $\OB$ planes, seven share the same index not shared by the remaining $\OB$ plane because each index occurs seven times. And since $\Phi_P$ is the dual of $\Phi_O$ given by multiplier of the remaining $\OB$ plus this shared index then any rotation from $\Rot{Q}$ will swap $\OB$ and $\PB4$ planes. The claim is that this can happen and the example shows it occurs in half of the $\Rot{Q}$ rotations. It is a conjecture that it always occurs for $\SB$ calibrations for at least the set $\Rot{Q}$.
\end{proof}
The dual of $\Phi$, $\Phi\hod$, is a co-non-associative calibration, which is the 8-form
\eqb{dPhi}\begin{split}
  \Phi\hod = &-\e{123456789ABCDEF}\Phi \\
     =&\quad \e{12478BDE} +\e{12479ACF} +\e{12568BCF} +\e{12569ADE} +\e{13468ADF}  \\
     &+\e{13469BCE} +\e{13578ACE} +\e{13579BDF} +\e{234589EF} +\e{2345ABCD} \\
     &+\e{236789CD} +\e{2367ABEF} +\e{456789AB} +\e{4567CDEF} +\e{89ABCDEF}\es
\end{split}\eqe


\begin{theorem}[Algebra stacking]\label{thm:2}
Each Cayley-Dickson algebra of $N$ generators has $H_N$ quaternion rings and is composed of $2^N-1$ embedded subalgebras with $N-1$ generators that each consist of $H_{N-1}$ embedded quaternion rings, where $H_N=(2^N-1)(2^N-2)/6$.
\end{theorem}
\begin{proof}
Each Cayley-Dickson algebra of $N$ generators embeds $H_N$ quaternion rings~\cite{Wilmot2} and $T_N=(2^N-1)(2^N-2)(2^N-4)/168$ Fano plane algebras~\cite{Wilmot2}. Each Fano plane algebra has seven quaternion rings, hence the number of quaternions shared between each Fano plane algebra is
  \[\frac{7T_N}{H_N} = \frac{6(2^N-4)}{24} = 2^{N-2}-1\ec\quad\text{for $N>2$.}\]
This is three for sedenions so that the quaternions in $\Theta$ match the quaternions in $\Phi$, $35=15\times7/3$. The quaternion rings shared amongst Fano plane subalgebras must also be shared amongst other larger subalgebras that cover the algebra so that dividing the shared quaternion rings in the algebra of $N$ generators with the quaternion rings for $N-1$ generators gives the number of embeddings of the largest subalgebra for $N>1$,
\eab{ll} H_N(2^{N-2}-1)/H_{N-1}&=(2^N-1)(2^N-2)(2^{N-2}-1)/(2^{N-1}-1)(2^{N-1}-2)\\
   &=(2^N-1)(2^N-2)(2^N-4)/(2^N-2)(2^N-4)\\
   &=2^N-1\es
\eae
The basis of the inductive argument is that, for $N=2$, $\HB$ embeds $\C$ $2^2-1=3$ times. For $N=3$, $\OB$ embeds $\HB$ $2^3-1=7$ times. It follows that $\SB$ embeds $2^4-1=15$ Fano-plane algebras. Each subalgebra is defined by a generalised 3-form calibration of $H_{N-1}$ terms that map to this many quaternions rings, defining a Cayley-Dickson algebra or quasi-algebra with $N-1$ generators.
\end{proof}
The embedded algebras can be found by first constructing the first quaternion calibration in $\Cl{2^N-1}$, $N>2$ by hand. It is easiest to extend the calibration 3-form in $\Cl{2^{N-1}-1}$, $N>3$ using the rule that all pairs of indices can only appear once. Then taking odd powers of this calibration until $(2^{N-1}-1)$-forms appear and selecting the subset of these that only share quaternions $2^{N-1}-1$ times. These $(2^{N-1}-1)$-forms extend $(2^{N-2}-1)$-forms of the calibration of $\Cl{2^{N-1}-1}$. If these calibrations, like the 8-form calibration in $\Cl{15}$, form a projection operator then this can only happen if the terms of its dual form a closed algebra under multiplication and the following definition is justified.

\begin{definition}[Sharp and Harp algebras]
The terms of the dual of the $(2^{N-1}-1)$-form calibration in $\Cl{2^N-1}$, along with unity, 1, form a subalgebra of $\Cp{2^N-1}$ denoted $\Sharp{N}$ (ignoring the calibration representation). The terms of the calibration and the imaginary pseudoscalar in $\Cl{2^N-1}$, along with $\Sharp{N}$ form a subalgebra of $\Cl{2^N-1}$ denoted $\Harp{N}$. $\Sharp{N}$ consists of even forms and $\Harp{N}$ consists of both even and odd forms.
\end{definition}
\begin{theorem}\label{thm:3}
The terms of $\Sharp{N}$ for $N\le4$ act under conjugation on itself to swap $2^{N-1}-1$ pairs of terms of an $(2^{N-1}-1)$-form calibration. They have the same action on $\Harp{N}$ and the associated 3-form calibration with $H_N$ terms. Hence they are related to automorphisms of $\Sharp{N}$ and the associated Cayley-Dickson algebra.
\end{theorem}
\begin{proof}
This theorem is trivial for $N=1,2$ and has been validated in~\cite{Wilmot1} for $N=3$ ($\OB$). The octonion case considers the action of $\Shar{i}3$ on $\theta_i$. As shown in \eqn{cal7} and demonstrated by the Fano plane, \fig{6sim}, each term of $\theta_i$ shares a basis vector with two other terms so multiplying $\theta_i\hod$ by any term $\theta_i\hod$ swaps three pairs of terms. Since the terms of $\GT$ are derived from these as three cycles of pairs of rotations then the action of the automorphisms are a subset of two of these pair swaps, cycling through the three possibilities generating the 21 terms of $\GT$. 

For sedenions, multiplying $\Phi\hod$ by a term of $\Phi_i\hod$ will swap seven pairs of 8-forms in $\Phi\hod$, which is summarised in \tar{inv}. From \eqn{Phi}, $\Phi_O$ terms contain $\e8$ whereas $\Phi_A$ and $\Phi_P$ terms do not. So multiplying by $\Phi_A\hod=\e{89ABCDEF}$ swaps all of $\Phi_O$ and $\Phi_P$. Multiplying by a member of $\Phi_P\hod$, $\Phi_i\hod$, $i\in\N_9^{15}$, will also swap between $\Phi_O$ and $\Phi_P$ apart from $\Phi_{i}$, which becomes the pseudoscalar while its partner, given by \eqn{obc1}, becomes $\e{1234567}$. The analogous situation occurs when multiplying by $\Phi_i\hod$, $i\in\N_2^8$, where a member of $\Phi_P$ swaps with $\Phi_A$ but all other swaps remain within their domain.

\tab{[ht]}{inv}{Pair swaps of $\Phi_i\hod$ acting on $\Phi$}
{\vskip 1pt\renewcommand{\arraystretch}{1.1}
\begin{tabular}{|c|c|c|c|c|} \hline
   $\times\Phi\hod$ terms &$\Phi_A\leftrightarrow$ &Within $\Phi_O$ &Within $\Phi_P$ &$\Phi_O\leftrightarrow\Phi_P$ \\\hline
  $\Phi_A\hod$ &0 &0 &0 &7 \\\hline
  $\Phi_O\hod$ &$\Phi_P$ &3 &3 &0 \\\hline
  $\Phi_P\hod$ &$\Phi_O$ &0 &0 &6 \\\hline
\end{tabular}}\tae
\end{proof}
Reference~\cite{Wilmot1} showed that combinations of rotations of the terms of the octonion dual calibration generated $\GT$ because such combinations are subset of the action of $\Sharp3$. The next section expands this concept to combinations of the action of $\Sharp4$ leading to all automorphisms of sedenions. The application of this theorem for cases $N>4$ is left as a conjecture. The action of $\Sharp{N}$ on itself is obvious so it only remains to show existence of the $(2^{N-1}-1)$-form calibration and is relation to the associative calibrations that map to Cayley-Dickson algebras. The combinations of 2-forms will always be a subset of this action thus defining automorphisms.

\section{Automorphisms of the Sharp algebras}
The $\Spin{n}$ group has basis rotations $\cos(\phi) +\e{ij}\sin(\phi)$ and the $\Pin{n}$ extends this to include reflections as combinations of basis elements $\e{i}$, for $i,j\in\N_1^n$ under the conjugation operation,
    \[\phi_x(y) = x y x^{-1}\text{, for $x,y\in\Pin{n}$.}\]
Automorphisms form a group and satisfy
\eqb{aut} \begin{split}
  \phi(xy) &= \phi(x)\phi(y)\text{ \quad and} \\
  \phi\circ\theta(x) &= \phi(\theta(x))\text{ \quad(composition),} 
\end{split} \eqe
and it is obvious that $\Pin{n}$ satisfies the first property due to transformations being conjugations when acting on $\Cl{n}$. Quaternions in $\Sharp2$ are their own automorphisms and this extends to $\Cl3=\Harp2$.

Integer fractions of $2\pi$ rotation angles form subalgebras of $\Spin{n}$ that are inner automorphisms of polyhedrons and certain reflections are outer automorphisms, together isomorphic to the Dihedral group. The members of the Sharp algebras consists of $180\degs$ rotations. 
A restriction to $90\degs$ conjugation rotations gives
\eqb{r90} \Rot{\e{ij}} = \cos\Big(\frac{\pi}4\Big) +\e{ij}\sin\Big(\frac{\pi}4\Big)\ec\eqe
and a pair of orthogonal such rotations can be divided into a 0-form and 4-form part plus 2-form parts, as
  \[\Rot{\e{ij}+\e{kl}}=\Rot{\e{ij}}\Rot{\e{kl}} = \alpha +\beta\ec\]
where $\alpha=\frac12(1+\e{ijkl})$ and $\beta=\frac12(\e{ij}+\e{kl})$.
These are shown in~\cite{Wilmot1} to have the rules $\alpha^2=\alpha$, $\beta^2=\alpha-1$ and $\alpha\beta=\beta\alpha=0$, so the $\Spin7$ rotation is then
\elb{rot7} 
  \Rot{\beta}x\Rot{\beta}^{-1}=&(\alpha+\beta)x(\alpha-\beta) \\
       =& \alpha x\alpha -\beta x\beta +\beta x\alpha -\alpha x\beta\es
\ele
For $x$ commuting with $\alpha$ and substituting $\alpha=1+\beta^2$ then
\eqb{roc7} \Rot{\beta}x\Rot{-\beta}=\beta(\beta x-x\beta) +x\es \eqe
This is the case for $\alpha\in\Shar{i}3$ and then $\beta\in\GT_i$ so that the calibration for this representation of $\GT$, $\theta_i$, satisfies $\theta_i\beta=\beta\theta_i$ and $\Rot{\beta}{\theta_i}=\theta_i$. Hence $\GT$ are automorphisms of $\Harp3$ and since \eqn{map7} maps the calibration into octonions, $\OB$, then $\GT$ are the automorphisms of $\OB$.

Equation \eqn{inv7} was shown to define an idempotent in~\cite{Wilmot1} and for the representation in \eqn{cal7} we have 
\[ \Theta_1^+=(1+\theta_1\hod)/8 = (1+\e{1256})(1+\e{1346})(1+\e{4567})/8\es \]
Most other terms of $\theta_1\hod$ can be used in the expansion but not all and different representations will be found including by changing the signs above. A form such as $\Theta_1^+$ is specified to be an ideal~\cite{Lounesto} and has the Radon-Hurwitz number $r_{-7}=r_1-4=-3$ so the order of the group is $2^3$ and the three terms in the ideal mean that $\Theta_1^+$ is a primitive ideal.

This is now generalised and extended to the complementary idempotent of $\Har{i}3$, $i\in\N_1^{30}$, for $k$ being the parity of $\theta_i$, which is the number of minus signs, as 
\elb{ideal7}
\Theta_i^+&=(1+(-1)^k\theta_i\hod)/8 \text{ and}\\
\Theta_i^-&=(7-(-1)^k\theta_i\hod)/8 = 1 - \Theta_i^+\ec
\ele
where $(\Theta_i^{\pm})^2=\Theta_i^{\pm}$ and $\Theta^+_i\Theta^-_i = 0$. 

Since any term with sign from $\theta_i\hod$ swaps terms of $\theta_i\hod$ under multiplication, apart from itself, which becomes 1, then $\Theta_i^+x=\Theta_i^+$ for $x$ a term of $\theta_i\hod$. This can be extended to $(1+x)/2$ that are factors in the expansion of $\Theta_i^+$ and since $\beta^2=\alpha-1$ then $2\beta^2-1$ gives these factors for $\beta\in\GT$. Hence
\[ \Theta_i^+(2\beta^2 +1)=\Theta_i^+\text{, for $\beta\in\GT_i$,}\]
and these factors form a basis and quotients of the $\Theta_i^+$ ideal. We can do more than find $\beta$ commutes with $\theta_i\hod$. For any $\beta\in\GT_i$, one of the terms of $\theta_i\hod$ is $2\alpha-1$ and $\alpha\beta=0$ so $\theta_i\hod\beta=-\beta+\delta$. In the orthogonal matrix representation of $\theta_i\hod$, where each term is four reflections and each reflected dimension occurs four times leaving three unreflected dimensions, then $\theta_i\hod$ has orthogonal matrix, -I, the negated identity matrix. Since one of the terms when multiplied by $\beta$ results in $\beta$ then the other terms must cancel so that $\delta=0$, giving
\[ \Theta_i^+\beta = 0\text{\quad and\quad}\Theta_i^-\beta = \beta\es \]
This is obvious since the ideals are complementary, $\Theta_i^++\Theta_i^-=1$, but only the 2-forms of $\GT$ cancel for $\Theta_i^+$, which is another way of defining $\GT$ as the automorphisms of $\Sharp3$. The $\Theta_i^-$ idempotent forms a null ideal space so $\Cl7$ cannot be decomposed into independent components.

The $\GT$ algebra as pairs of rotations transforms $\theta_i$ to itself, meaning that the labels of the 6-simplex are invariants. This is not the case for the Dihedral group whereby reflections change the parity of labelled polyhedra. The two are distinguished as inner and outer automorphisms where outer automorphisms that change the labelling cannot be performed as conjugation operations. In the Pin group both transformations use conjugation so the terminology is not appropriate. New terminology is introduced to cover conjugation reflections, some of which will be found to be automorphisms, but cannot be called outer automorphisms.
\begin{definition}\label{def:weak}
Any isomorphism from $\Cl{n}$ that maps the labelled, independent simplex faces to another simplex and keeps the labels of the simplex in the same relative order shall be called a {\it strong isomorphism}, which implies the associative calibration is invariant. Any transformation that changes the calibration to another calibration is just a change of labels so is an isomorphism herein called a {\it weak isomorphism}. This is divided into {\it internal} isomorphisms that exist in the same primary and {\it external} whereby the labels change to equivalent labels in another primary. Otherwise the conjugate transformation will produce a quasi-calibration and will not be an isomorphism. All isomorphisms derived from $\Phi$ will be denoted with superscript indices to distinguish them from the $\Phi=\Phi_{A+O+P}$ and $\Phi\hod$ calibration domains.
\end{definition}
All single reflections from the respective Pin algebra will change the parity of some independent faces or signs in the associative calibration. Hence these are weak isomorphisms as are the pseudo-rotations and pseudoscalar itself will be found later to be both an automorphism and a weak, internal isomorphism. It is obvious that strong isomorphisms are automorphisms but this is not a bijection.

In $\Cl3$, the $90\degs$ rotations from $\Spin7$ are weak, external isomorphisms. The members of $\GT$ provide another algebra element that restores the original calibration to define automorphisms. But the opposite parity of these members with the second rotation negated will also return the calibration to the same primary and so are weak, internal isomorphisms. Hence $\GT$ can be extended by adding $\GT$ members with opposite parity, denoted $\GT'$, which generates the algebra $\GT+\GT'$ with Lie product summarised in \tar{g22}.
\tab{[h!]}{g22}{$\GT$ and $\GT'$ Lie product summary}
\begin{tabular}{|c|c c|} \hline
  [.,.]      &$\GT$      &$\GT'$ \\\hline
 $\GT$  &$\GT$      &$\GT+\GT'$  \\
 $\GT'$ &$\GT+\GT'$ &$\GT$ \\\hline
\end{tabular} \tae

\thr{3} showed that terms of $\Shar13$ swap terms of the calibration given in \eqn{cal7}, which means they commute with the calibration and hence are strong isomorphisms of $\Har13$. The calibration $\theta_1$ is an invertible form, \eqn{inv7}, so that all members of $\Pin{7}$ under conjugation keep this property invariant and thus define calibrations, but these are weak isomorphisms of $\Har13$. To be automorphisms they have to form a group and satisfy \eqn{aut}. The reflections are isomorphic to $\ZB{7}\cong\Z/7\Z$ and the pseudoscalar extends this to the cycle of sixteen internal calibrations. But only the identity and the pseudoscalar satisfy \eqn{aut}. The 21 elements of $\GT$ were derived in~\cite{Wilmot1} from what was called the effective algebra, now called $\Shar13$ due to \thr{3}, as pairs of 2-forms. The three ordered pairs of the indices of any 4-form term of $\theta_1\hod$ generate three rotations each, which act as a subset of the invariant action of $\Shar13$ on $\theta_1$, so are also invariants. Cycling through all three rotations as combinations of 2-form pairs in any term of $\theta_1\hod$ swaps all three combinations of 3-form pairs having the same effect as 4-form conjugation swapping three pairs of 3-forms. These 21 pairs of rotations are the basis of $\GT$, which has 14 independent basis elements shown in \tar{gt7} for the representation of $\e{3568BDE}$. This can be extended to $\GT\times\C$ where the imaginary is represented by the pseudoscalar but $\GT\times\e{1234567}$ are weak isomorphisms.
\tab{[ht]}{gt7}{$\GT$ representation for $\e{3568BDE}$ from $\Phi_O$}
\begin{tabular}{cclcclccl}
  A&=&$\frac12(\e{56} -\e{8B})$ &B&=&$\frac12(-\e{36} -\e{8D})$ &C&=&$\frac12(\e{35} +\e{8E})$ \\[2pt]
  D&=&$\frac12(-\e{3B} +\e{5B})$ &E&=&$\frac12(\e{38} -\e{5E})$ &F&=&$\frac12(\e{3E} +\e{58})$ \\[2pt]
  G&=&$\frac12(-\e{3D} -\e{5B})$ &H&=&$\frac12(\e{8B} -\e{DE})$ &I&=&$\frac12(\e{8D} +\e{BE})$ \\[2pt]
  J&=&$\frac12(-\e{8E} +\e{BD})$ &K&=&$\frac12(-\e{5D} -\e{6E})$ &L&=&$\frac12(-\e{5E} +\e{6D})$ \\[2pt]
  M&=&$\frac12(-\e{3E} +\e{6B})$ &N&=&$\frac12(\e{5B} -\e{68})$ &&&
  \end{tabular}\tae

Double $90\degs$ $\Sharp3$ rotations were shown in~\cite{Wilmot1} to generate $\GT$ using \eqn{roc7}. 
Here we state the equivalent results for $\Sharp{4}$ by extending $\Rot\beta$ in \eqn{r90} to quadruple $90\degs$ rotations, $\beta=(\e{ij}+\e{kl}+\e{mn}+\e{op})/2$,
\[ \begin{split}
\Rot{\beta} 
  &= (1+\e{ij})(1+\e{kl})(1+\e{mn})(1+\e{op})/4 \\
  &= (1+\e{ij}+\e{kl}+\e{mn}+\e{op}+\e{ijkl}+\e{ijmn}+\e{ijop}+\e{klmn}+\e{klop}+\e{mnop}\\
  &\qquad+\e{ijklmn}+\e{ijklop}+\e{ijmnop}+\e{klmnop}+\e{ijklmnop})/4\es \\
\end{split}\]

For $\Spin{15}$ the full rotation is divided up into 0 plus 8-form, 2-form, 6-form and 4-form parts as half of $\alpha, \beta,\gamma, \delta$, respectively. The full rotation is, $x'=\Rot{\beta}x\Rot{-\beta}$, 
\elb{rot15} x^\prime = (\alpha+\beta+\gamma+\delta)x(\alpha-\beta-\gamma+\delta)/4\ec \ele
where 
\begin{center}\vskip -4ex\begin{tabular}{ll} $\beta^2 = \delta -1\ec$ &$\alpha^2=\alpha\ec$ \\
    $\delta^2 = 3\alpha -2\delta\ec$ &$(\alpha+\beta)^2 = \alpha +\beta +\delta -\gamma -1\ec$ \\
    $(\beta+\gamma)^2 = 4\alpha -4\ec$ &$(\alpha+\delta)^2=4\alpha\es$
\end{tabular} \end{center}
Defining $\alpha^\prime=\alpha+\delta$ and $\beta^\prime=\beta+\gamma$ with properties $\alpha^{\prime2}=4(\alpha^\prime-\delta)$, $\beta^{\prime2}=4(\alpha^\prime-\delta-1)$ and $\alpha^\prime\beta^\prime= \beta^\prime\alpha^\prime=0$ then the rotation becomes, for $x\in\Harp{4}$,
\[ \begin{split} x^\prime&=(\alpha^\prime+\beta^\prime)x(\alpha^\prime-\beta^\prime)/4 \\
  &= (\alpha^\prime x\alpha^\prime -\beta^\prime x\beta^\prime +\beta^\prime x\alpha^\prime -\alpha^\prime x\beta^\prime)/4\es
\end{split} \]
For $[\alpha^\prime,x]=0$ and substituting $\alpha^\prime=\beta^{\prime2}/4 +\delta +1$, then
\[\begin{split}
  x^\prime&=(\alpha^\prime-\delta)x -\beta^\prime x\beta^\prime/4 \\
    &= (\beta^{\prime2}/4+1)x-\beta` x\beta`/4 \\
    &= \beta^\prime(\beta^\prime x-x\beta^\prime)/4 +x\es
\end{split} \]
Similar to 7-D this is also expected to be satisfied by the calibration. We find both $\alpha$ and $\delta$ commute with all terms of $\Harp{4}$~\cite{Wilmot3} so $\alpha^\prime$ terms are strong isomorphisms. 
It is instructive to apply the same automorphism analysis from $\Sharp3$ to $\Sharp{4}$ in order to study the ideals and later apply the results to sedenions. Like $\Sharp3$ in \eqn{inv7}, the calibration needs to be negated to get negative parity so only positive parity of $\Phi$ is used for the ideals.

The dual of $\psi$ from \thr{1} means $(7-\Phi\hod)^2=64$ so that the ideals for $\Phi$ are
\elb{id15}
\Psi^+ &=(1+\Phi\hod)/16 = (1+\e{12478BDE})(1 +\e{12479ACF}) \\
       & \qquad\qquad\qquad\qquad(1+\e{12569ADE})(1 +\e{13469BCE})/16\text{ and}\\
\Psi^- &=(15-\Phi\hod)/16 = 1 -\Psi^+\\
\ele
where $(\Psi^{\pm})^2=\Psi^{\pm}$ and $\Psi^+\Psi^- = 0$. Hence $\beta^\prime$ that satisfies $\Psi^-\beta^\prime=\beta^\prime$ will be a strong isomorphisms of $\Sharp{4}$, $\Phi^\prime=\Phi$. The members of the ideal are the Harp algebra,
\eab{cl} \Psi^+x &= \Psi^+\text{ and} \\
     \Psi^+\e{123456789ABCDEF}x &= \e{123456789ABCDEF}\Psi^+ \text{, for all $x\in\Sharp{4}$.}
\eae
Here, four terms of $\Phi\hod$ have been turned into sub-idempotents, such as $\frac12(1+\e{12478BDE})$, and four idempotents multiplied to recover $\Psi^+$. Almost any four terms of $\Phi\hod$ constructed in this way generates the ideal. Of the $\binom{15}4=1365$ combinations, 840 satisfy this condition and 525 generate other calibrations.

Similar to $\Theta_1^+$ in $\Shar{i}3$, $\Psi^-$ is specified to be an ideal~\cite{Lounesto} and has the Radon-Hurwitz number $r_{-15}=r_{-7}-4=r_1-8=-7$ so the order of the group is $2^7$ and hence $\Psi^+$ is not a primitive ideal. 

Four $90\degs$ rotations correspond to a product of eight vectors, which is a Pfaffian of degree~7 containing 105 terms. Expanding the last column gives seven terms involving the last vector with products of each other vector containing a Pfaffian of degree~5 containing 15 terms. In order to keep $\Phi$ and $\Theta$ within the same primary representation, these can be decreased to seven cases.

\begin{theorem}\label{thm:4}
The invariants of $\Phi$ under $90\degs$ rotations are given by two maps of each term of $\Phi\hod$ in \eqn{dPhi}. The first, called cyclic because it is based upon the $\GT$ map for each half of the 8-forms, is
\eqb{aut1}\begin{split}
\e{\mu_1\mu_2\mu_3\mu_4\mu_5\mu_6\mu_7\mu_8}\longrightarrow
&(\e{\mu_1\mu_2}+\e{\mu_3\mu_4}+\e{\mu_5\mu_6}+\e{\mu_7\mu_8},\\
&\ \e{\mu_1\mu_3}+\e{\mu_2\mu_4}+\e{\mu_5\mu_7}+\e{\mu_6\mu_8}, \\
&\ \e{\mu_1\mu_4}+\e{\mu_2\mu_3}+\e{\mu_5\mu_8}+\e{\mu_6\mu_7})\es
\end{split}\eqe
where $\mu_i$ is the $i$\textsuperscript{th} index of $\e{\mu_1\mu_2\mu_3\mu_4\mu_5\mu_6\mu_7\mu_8}$. This map avoids the quasi-cal\-ibrations provided in the corollary to \thr{1}. Also note that $\Phi_O$ corresponds to this structure and other combinations that mix the lower or upper indices will change the primary representation. The same argument applies for $\Phi_P$ and only the lower indices apply to $\Phi_A$, which limits this map to four cases. The second map, is called mixed because it mixes some 7-D Fano plane indices with false Fano plane indices, and $\e8$ for $\Phi_O$, which is
\eqb{aut2}\begin{split}
\e{\mu_1\mu_2\mu_3\mu_4\mu_5\mu_6\mu_7\mu_8}\longrightarrow\ 
    &(\e{\mu_1\mu_5}+\e{\mu_2\mu_6}+\e{\mu_3\mu_7}+\e{\mu_4\mu_8}, \\
    &\ \e{\mu_1\mu_6}+\e{\mu_2\mu_5}+\e{\mu_3\mu_8}+\e{\mu_4\mu_7}, \\
    &\ \e{\mu_1\mu_7}+\e{\mu_2\mu_8}+\e{\mu_3\mu_5}+\e{\mu_4\mu_6}, \\
    &\ \e{\mu_1\mu_8}+\e{\mu_2\mu_7}+\e{\mu_3\mu_6}+\e{\mu_4\mu_5})\es
\end{split}\eqe
This is necessary in order to consider $\Phi_A$ swapping with other terms in $\Phi$. This may introduce quasi-calibrations that will be excluded later. For the last term in the Pfaffian to be $\e{\mu_i\mu_8}$, $i\in\N_1^4$, then there is only one case with all lower and upper indices, thus giving four cases from the Pfaffians of degree~5, shown in \eqn{aut2}. Applying these seven maps to the fifteen terms of $\Phi\hod$ gives 105 first sign variation invariants, which are provided in the appendix in \tas{auta}{aut2p}. We now consider paired negation of terms, since each Pfaffian vector may be a $\pm90\degs$ rotation and since the calibration parity needs to be preserved there are only four cases,
\begin{align*}
\e{\mu_1\mu_2}+\e{\mu_3\mu_4}+\e{\mu_5\mu_6}+\e{\mu_7\mu_8}\longrightarrow\ 
    &(\e{\mu_1\mu_2}+\e{\mu_3\mu_4}+\e{\mu_5\mu_6}+\e{\mu_7\mu_8}, \\
    &\ \e{\mu_1\mu_2}+\e{\mu_3\mu_4}-\e{\mu_5\mu_6}-\e{\mu_7\mu_8}, \\
    &\ \e{\mu_1\mu_2}-\e{\mu_3\mu_4}+\e{\mu_5\mu_6}-\e{\mu_7\mu_8}, \\
    &\ \e{\mu_1\mu_2}-\e{\mu_3\mu_4}-\e{\mu_5\mu_6}+\e{\mu_7\mu_8})\ec
\end{align*}
quadrupling the number of sign variations of each table to generate 420 isomorphisms of $\Phi$ and $\Theta$.
\end{theorem}

\begin{proof}
The first cyclic map, \eqn{aut1}, is identical to the one proved by Schafer~\cite{Schafer} for $\Phi_O$ as automorphisms of sedenions, discussed later. This was extended with the first sign variation by Brown~\cite{Brown} whereas here the new 7-form calibration, $\Phi=\Phi_A +\Phi_O +\Phi_P$, is considered first. The second map, \eqn{aut2}, extends the first map to cover all groups of two of the eight indices of each $\Sharp{4}$ term (apart from the scalar). But the constraint is we are considering terms of $\Phi\hod$ in $\eqn{dPhi}$ that map to other terms. Apart from $\Phi_A\hod=\e{89ABCDEF}$ this involves keeping the 7-D indices and $\e{8}$ plus false Fano plane indices separate. Hence, once $\e{\mu_1\mu_i}$, $i\in\N_2^8$ is selected from $\e{\mu_1\mu_2\mu_3\mu_4\mu_5\mu_6\mu_7\mu_8}$, then for $i\le4$ the second term must also select from indices $\mu_3$ and $\mu_4$ in $\N_2^4$. The same argument applies to the upper four indices. To cover $\Phi_A\hod$, if $i>4$ then the second term must have $\mu_3\in\N_2^4$ and $\mu_4\in\N_5^8$, etc. This means that the first term, $\e{\mu_1\mu_i}$ for $i\in\N_2^8$ governs the solution and there are only seven first sign variation invariants for each term in $\Phi\hod$, giving a total of 105 before the sign variations are applied. Unlike $\Sharp3$, the parity can cancel under the map \eqn{map15} so all even permutations of the four $90\degs$ rotation 2-forms should be considered. Half of these permutations will be negations of the other half so only four need be considered. 

\tar{inv} shows the action of $\Sharp{4}$ rotations on $\Phi$ as swaps of seven pairs of terms. The four $90\degs$ rotations of each $\Phi^A$ terms are a subset with swaps of four pairs of terms, as shown in \tar{aut}.
Thus there are 420 elements of $\Sharp{4}$ that keep $\Harp{4}$ invariant, so are automorphisms of $\Harp{4}$ if they form a group and satisfy \eqn{aut}.
\end{proof}
The Mixed($\o4$) row of \tar{aut} denotes the $\Phi_O^M$ invariant terms that map to $\o4$. Since $\o4$ is mapped from $\e8$ in \eqn{map15} then such transformations do not rotate octonion components within $\Phi_O$. The following Mixed(other) row in the table has terms that do not map to $\o4$ so acts on $\Phi$ differently. 
\tab{[ht]}{gt15}{$\Cl{15}$ representation of $\GT$}
\begin{tabular}{cclccl}
  $\GG{A}1$&=&$\frac12(\e{23}-\e{45}+\e{AB}-\e{CD})$
  &$\GG{B}1$&=&$\frac12(-\e{13}-\e{46}-\e{9B}-\e{CE})$\\[2pt]
  $\GG{C}1$&=&$\frac12(\e{12} -\e{47} +\e{9A} -\e{CF})$
  &$\GG{D}1$&=&$\frac12(-\e{15}+\e{26}-\e{9D}+\e{AE})$\\[2pt]
  $\GG{E}1$&=&$\frac12(\e{14}+\e{27}+\e{9C}+\e{AF})$ 
  &$\GG{F}1$&=&$\frac12(-\e{17}+\e{24}-\e{9F}+\e{AC})$\\[2pt]
  $\GG{G}1$&=&$\frac12(-\e{16}-\e{25}-\e{9E}-\e{AD})$
  &$\GG{H}1$&=&$\frac12(\e{45}+\e{67}+\e{CD}+\e{EF})$\\[2pt]
  $\GG{I}1$&=&$\frac12(\e{46}-\e{57}+\e{CE}-\e{DF})$
  &$\GG{J}1$&=&$\frac12(\e{47}+\e{56}+\e{CF}+\e{DE})$\\[2pt]
  $\GG{K}1$&=&$\frac12(-\e{26}+\e{37}-\e{AE}+\e{BF})$ 
  &$\GG{L}1$&=&$\frac12(\e{27}+\e{36}+\e{AF}+\e{BE})$\\[2pt]
  $\GG{M}1$&=&$\frac12(\e{17}+\e{35}+\e{9F}+\e{BD})$
  &$\GG{N}1$&=&$\frac12(\e{25}-\e{34}+\e{AD}-\e{BC})$ 
\end{tabular}\tae
\begin{definition}
\thr{4} finds three cyclic \eqn{aut1} and four mixed \eqn{aut2} isomorphisms of $\Harp{4}$. These are identified as $\Phi^C$ and $\Phi^M$, respectively, and together denoted $\Phi^A$. 
The theorem also finds four sign variations as $(+,+,+,+)$, $(+,+,-,-)$, $(+,-,+,-)$, $(+,-,-,+)$, which are denoted $\Phi^{A(1)}$, $\Phi^{A(2)}$, $\Phi^{A(3)}$, $\Phi^{A(4)}$, respectively. Hence the total signed transformations of $\Harp{4}$ are denoted as $\Phi^A=\Phi^{A(1+2+3+4)}$.
\end{definition}
The Lie product tables for the invariants $\Phi_A^A$, $\Phi_O^C$, $\Phi_P^C$, $\Phi_O^M$ and $\Phi_P^M$ are considered separately and these will be called isomorphism domains. They require all four sign variations to check for closure under the Lie product, where closure means the Lie product of terms from $\Phi^A$ result in terms from $\Phi^{A(1+2+3+4)}$.
The first sign variation of these domains are \tas{auta}{aut2p}, respectively, in the appendix, and the sign variations have been juggled to provide signed domain closure and cross domain closure, as much as possible.

\tab{[!ht]}{aut}{Pair swaps of $\Phi$ transformations acting on $\Phi$}%
{\renewcommand{\arraystretch}{1.1}%
\begin{tabular}{|cc|c|c|c|c|} \hline
   $\times\Phi\hod$ terms &\multicolumn{1}{|c|}{$\times\Phi^A$ terms}&$\Phi_A\leftrightarrow$ &Within $\Phi_O$ &Within $\Phi_P$ &$\Phi_O\leftrightarrow\Phi_P$ \\\hline
  \multirow{3}{*}{$\Phi_A\hod$} & &\bf0 &\bf0 &\bf0 &\bf7 \\\cline{2-6}
  &\multicolumn{1}{|c|}{Cyclic} &0 &0 &0 &4\\\cline{2-6}
  &\multicolumn{1}{|c|}{Mixed} &0 &0 &0 &4\\\hline
  \multirow{3}{*}{$\Phi_O\hod$} & &$\Phi_P$ &\bf3 &\bf3 &\bf0 \\\cline{2-6}
  &\multicolumn{1}{|c|}{Cyclic} &0 &2 &2 &0\\\cline{2-6}
  &\multicolumn{1}{|c|}{Mixed($\o4$)} &1 &0 &3 &0\\\cline{2-6}
  &\multicolumn{1}{|c|}{Mixed(other)}  &1 &2 &1 &0\\\hline
  \multirow{3}{*}{$\Phi_P\hod$} & &$\Phi_O$ &\bf0 &\bf0 &\bf6 \\\cline{2-6}
  &\multicolumn{1}{|c|}{Cyclic} &0 &0 &0 &4\\\cline{2-6}
  &\multicolumn{1}{|c|}{Mixed} &1 &0 &0 &3\\\hline
\end{tabular}}\tae
The first map, \eqn{aut1}, applied to $\Phi_O\hod$ was based on $\GT$ so it is appropriate to duplicate the signs to 15 dimensions, as shown in \tar{gt15}, which forms a 15-dimensional representation of $\GT$, since its construction is the same as that of Schafer~\cite{Schafer}. The 21 invariants of $\GT$, $\Phi_O^{C(1)}$, are shown in \tar{aut1o} and the 21 invariants of $\Phi_P^{C(1)}$ are provided in \tar{aut1p}. The first 14 elements of \tar{aut1o} are independent providing rank 14 and $\Phi_{O+P}^{C(1)}$ has rank 30. The second map, which in general mixes the two $\GT$ representations, has two lots of 28 weak isomorphisms, those of $\Phi_O^{M(1)}$ are shown in \tar{aut2o} and those from $\Phi_P^{M(1)}$ in \tar{aut2p}. The combined maps for $\Phi_A$, $\Phi_A^{A(1)}$, are shown in \tar{auta}. 
An example Lie product for the 15-D $\GT$ is the first and second members of $\Phi_O^{C(1)}$, taken from \tar{aut1o}, as 
\[ \GG{A}1 = \frac12(\e{23}-\e{45}+\e{AB}-\e{CD}) \text{\ and\ } \GG{B}1 =\frac12(-\e{13}-\e{46}-\e{9B}-\e{CE})\ec \]
and the Lie product of these two elements is $[\GG{A}1,\GG{B}1]=-\GG{C}1 -\GG{J}1$, where these are the third and tenth elements of the table. The second sign variation, $\Phi_O^{C(2)}$, reverses the sign of the last two terms, which corresponds to Brown's~\cite{Brown} extension of Schafer's results. It should be noted that this is interpreted in the literature as $120\degs$ rotations between the 7-D Fano plane and the false Fano plane plus $\e8$. Here this extension to $\GT$ is derived from $90\degs$ rotations and \tar{aut} shows that all Cyclic $\Phi_O$ isomorphisms swap Fano planes within the $\Phi_O$ and $\Phi_P$ set of planes. Such transformations do not rotate the elements of $\e{1234567}$ so this result disagrees with the usual interpretation. 
Lie products between domains can involve other members of other algebras. For example, the Lie product of the fourteenth member of $\Phi_P^{C(1)}$, from \tar{aut1p}, with ${\rm A}^{(3)}$ is
\[ 
  \Big[\frac12(\e{24} -\e{35} +\e{8E} -\e{9F}),\frac12(\e{23} +\e{45} +\e{AB} +\e{CD})\Big]=\e{25}+\e{34}\es
\]
This is twice the fourteenth member of $\GT$ with a change of parity for the second term. 
\begin{definition}
The exceptional Lie algebra $\GT$ in 7-D for the calibration \eqn{cal7} will be designated as $\Phi^{G(1)}$ as provided in \tar{gt15}. This matches the first two terms of each element in \tar{aut1o}. Changing the parity of the second term of each pair is denoted $\Phi^{G(2)}$, which equals $\GT'$ from \tar{g22}, and extends $\GT$ to $\Phi^{G(1+2)}=\Phi^G$, with $\Phi^{G(1)}$ provided in the appendix, \tar{autg}.
\end{definition}
The summary Lie product table of these algebras from \tar{g22} is
\[ \left[\Phi^{G(i)},\Phi^{G(j)}\right] \rightarrow \begin{cases}
    \Phi^{G(1)}\ec    &\text{for $i\in\N_1^2$, $j=i$, and} \\
    \Phi^{G(1+2)}\ec  &\text{otherwise.} 
\end{cases} \]
\begin{definition}
The algebra $\Phi^{N(1+2)}=\Phi^N$ is introduced composed from the 44 terms of $\Phi^{N(1)}$, shown in the appendix, \tar{autn}, and the changed parity of the second 2-form term of each element, denoted $\Phi^{N(2)}$.
\end{definition}
The Lie product summary is similar to $\Phi^G$,
\[ \left[\Phi^{N(i)},\Phi^{N(j)}\right] \rightarrow \begin{cases}
    \Phi^{N(1)}\ec    &\text{for $i\in\N_1^2$, $j=i$, and} \\
    \Phi^{N(1+2)}\ec  &\text{otherwise.} 
\end{cases} \]
The algebras $\Phi^G$ and $\Phi^N$ are found in the domain cross products in the following.
\begin{lemma}\label{lem:2}
    The sum of members of $\Phi_A^A$, $\Phi_O^C$, $\Phi_P^C$, $\Phi^G$ and $\Phi^N$ form a closed algebra.
\end{lemma}
\begin{proof}
The Lie products of sign variations of $\Phi_A^{A}$ with itself generate the  results shown in \tar{aa}.
\tab{[ht]}{aa}{Lie Products of $\Phi_A^{A(1+2+3+4)}$}\centering \begin{tabular}{|c|c|c|c|c|}\hline
 \multicolumn{5}{|l|}{$\tstrut\big[\Phi_A^{A(i)},\Phi_A^{A(j)}\big]$} \\[.2ex]\hline
 \diagbox[height=3.5ex]{$i$}{$\ j$} &1 &2 &3 &4 \\\hline
1 &$\tstrut\Phi_A^{A(1)}$   &$\Phi_A^{A(1+2)}$ &$\Phi_A^{A(1+3+4)}$  &$\Phi_A^{A(3+4)}$\\
2 &$\tstrut\Phi_A^{A(1+2)}$   &$\Phi_A^{A(1)}$ &$\Phi_A^{A(1+3+4)}$  &$\Phi_A^{A(3+4)}$\\
3 &$\tstrut\Phi_A^{A(1+3+4)}$   &$\Phi_A^{A(1+3+4)}$ &$\Phi_A^{A(1+4)}$  &$\Phi_A^{A(1+3+4)}$\\
4 &$\tstrut\Phi_A^{A(3+4)}$   &$\Phi_A^{A(3+4)}$ &$\Phi_A^{A(1+3+4)}$  &$\Phi_A^{A(1+3+4)}$\\\hline
\end{tabular}\tae

The mixed part of $\Phi_A^A$ is included, even though $\Phi_A^C$ is closed under Lie products within itself, because $[\Phi_A^{C(1)},\Phi_O^{C(1)}]\rightarrow \Phi_A^A$.

The Lie products between $\Phi^G$, $\Phi^N$, $\Phi_A^A$ and $\Phi_{O+P}^C$ give the following results,
\[ \begin{split}
  \big[\Phi^{G(1+2)}, \Phi^{N(1+2)}\big] &\rightarrow 0\ec \\
  \big[\Phi^{G(1+2)}, \Phi_A^{A(1+2+3+4)}\big] &\rightarrow 0\ec \\
  \big[\Phi^{N(1+2)}, \Phi_A^{A(1+2+3+4)}\big] &\rightarrow \Phi_A^{(1+3+4)} +\Phi^{N(1)}\ec \\
  \big[\Phi^{N(1+2)}, \Phi_{O+P}^{C(1+2+3+4)}\big] &\rightarrow \Phi^{N(1+2)}\ec \\
  \big[\Phi_A^{A{(1+2+3+4)}},\Phi_{O+P}^{C(j)}\big] &\rightarrow
   \Phi_A^{A(1+3+4)} +\Phi^{N(1+2)}\ec \\
  \big[\Phi^{G(1+2)}, \Phi_O^{C(j)}\big] &\rightarrow 
    \begin{cases}
        \Phi^{G(1)}\ec &\text{if $j\in\N_1^2$ and} \\
        \Phi^{G(1+2)}\ec &\text{otherwise, and} 
    \end{cases} \\
  \big[\Phi^{G(1+2)}, \Phi_P^{C(j)}\big] &\rightarrow 
    \begin{cases}
        \Phi^{G(1+2)}\ec &\text{if $j\in\N_1^2$ and} \\
        \Phi^{G(1)}\ec &\text{otherwise.} 
    \end{cases} \\
\end{split}\]
The Lie products of $\Phi_{O+P}^C$ are shown in \tar{op} with the superscript C and brackets suppressed. Hence $\Phi_A^A +\Phi_O^C +\Phi_P^C +\Phi^G +\Phi^N$ forms a closed algebra. Since $\Phi_O^{C(1)}$ is a 15-D representation of $\GT$ then \tar{op} demonstrates that $\Phi_O^{C(1+2)}\cong\GT\times\Sym2$.
Also $\Phi_P^{C(1+2)}\cong\Phi_{O+P}^{C(1)}\times\Sym2$.

\tab{[ht]}{op}{Lie products of $\Phi_{O+P}^{C(1+2+3+4)}$}
\renewcommand{\arraystretch}{1.6}\setlength\tabcolsep{3pt}\small
\begin{tabular}{|l|c|c|c|c|c|c|c|c|}\hline
$[.,.]\tstrut$ &$\Phi_O^{1}$ &$\Phi_O^{2}$ &$\Phi_O^{3}$  &$\Phi_O^{4}$ 
  &$\Phi_P^{1}$    &$\Phi_P^{2}$    &$\Phi_P^{3}$       &$\Phi_P^{4}$ \\\hline
$\Phi_O^{1}\tstrut$ &$\Phi_O^1$ &$\Phi_O^1$ &$\Phi_O^{1\ssmp3}$ &$\Phi_O^{1\ssmp4}$
  &$\Phi_{O\smp P}^{1}$ &$\Phi_{O\ssmp P}^{1}$ &$\Phi^G\smp\Phi_O^{1}\smp\Phi_P^{3\ssmp4}$
  &$\Phi^G\smp\Phi_O^{1}\smp\Phi_P^{3+4}$ \\
$\Phi_O^{2}\tstrut$  &.  &$\Phi_O^{1}$ &$\Phi_O^{1\ssmp4}$ &$\Phi_O^{1\ssmp3}$
  &$\Phi_{O\ssmp P}^{1}$ &$\Phi_{O\ssmp P}^{1}$
   &$\Phi^G\smp\Phi_O^{1}\smp\Phi_P^{3\ssmp4}$  &$\Phi^G\smp\Phi_O^{1}\smp\Phi_P^{3\ssmp4}$ \\
$\Phi_O^{3}\tstrut$      &.  &.  &$\Phi_O^{1\ssmp3}$ &$\Phi_O^{1\ssmp4}$ 
  &$\Phi^G\smp\Phi_{O\ssmp P}^{3\ssmp4}\smp\Phi_P^{1}$   &$\Phi^G\smp\Phi_{O\ssmp P}^{3\ssmp4}\smp\Phi_P^{1}$
   &$\Phi_{O\ssmp P}^{3\ssmp4}\smp\Phi_P^{1}$ &$\Phi_{O\ssmp P}^{3\ssmp4}\smp\Phi_P^{1}$ \\
$\Phi_O^{4}\tstrut$      &.  &.  &.   &$\Phi_O^{1\ssmp3}$
   &$\Phi^G\smp\Phi_{O\ssmp P}^{3\ssmp4}\smp\Phi_P^{1}$   &$\Phi^G\smp\Phi_{O\ssmp P}^{3\ssmp4}\smp\Phi_P^{1}$   &$\Phi_{O\ssmp P}^{3\ssmp4}\smp\Phi_P^{1}$ &$\Phi_{O\ssmp P}^{3\ssmp4}\smp\Phi_P^{1}$ \\
$\Phi_P^{1}\tstrut$      &.  &.  &.   &.   
   &$\Phi_{O\ssmp P}^{1}$ &$\Phi_{O\ssmp P}^{1}$ &$\Phi_{O\ssmp P}^{3\ssmp4}\smp\Phi_P^{1}$ &$\Phi_{O\ssmp P}^{3\ssmp4}\smp\Phi_P^{1}$ \\
$\Phi_P^{2}\tstrut$      &.  &.  &.   &.   
   &. &$\Phi_{O\ssmp P}^{1}$ &$\Phi_{O\ssmp P}^{3\ssmp4}\smp\Phi_P^{1}$ &$\Phi_{O\ssmp P}^{3\ssmp4}\smp\Phi_P^{1}$ \\
$\Phi_P^{3}\tstrut$      &.  &.  &.   &.
   &. &. &$\Phi_O^{1}\smp\Phi_P^{3\ssmp4}$   &$\Phi_O^{1}\smp\Phi_P^{3\ssmp4}$ \\
$\Phi_P^{4}\tstrut$      &.  &.  &.   &.
   &. &. &. &$\Phi_O^{1}\smp\Phi_P^{3\ssmp4}$ \\\hline
\end{tabular}\tae
\end{proof}

\begin{definition}
  Define the sign variations of subalgebras from \ler{2} as
  \[ \begin{split}
    \Phi^{S(1)} &= \Phi_A^{A(1+3+4)} +\Phi_O^{C(1+3+4)} +\Phi_P^{C(1+3+4)} +\Phi^{G(1+2)} +\Phi^{N(1+2)}\text{ and} \\
    \Phi^{S(2)} &= \Phi_A^{A(2)} +\Phi_O^{C(2)} +\Phi_P^{C(2)}\es 
\end{split} \]
\end{definition}
\begin{theorem}\label{thm:5}
   The weak internal isomorphisms of $\Harp{4}$ are $\Phi^A$ and are the basis of two distinct algebras under the Lie product.
\end{theorem}
\begin{proof}
\thr{4} and \ler{2} showed that $\Phi^{S(1+2)}$ forms a group. The 4-form components also satisfy the ideal condition, $\Psi^-\beta'=\beta'$ where $\beta'=\beta+\gamma$, $\beta\in\Phi_A^A$, and $\gamma=(2\beta^3+5\beta)/3$ is the 6-form part of the rotation, for $\beta\in\{\Phi_A^A+\Phi_{O+P}^{C}\}$. Furthermore, from the corollary to \thr{1}, the 2-form terms of $[\Phi_A^A,\Phi_{O+P}^{C}]$ and $[\Phi_{O}^{C},\Phi_{P}^{C}]$ do not contain terms that transform the calibration to a quasi-calibration and being derived from terms of the dual calibration that satisfy $\Psi^-\beta'=\beta'$, are strong isomorphisms. The domain cross products, $\Phi^G$ and $\Phi^N$ do not keep \eqn{obc} invariant, which generates external calibrations but are derived from the internal isomorphisms so $\Phi_A^A +\Phi_O^C +\Phi_P^C$ is the basis of a group.

This leaves the mixed domains, $\Phi_{O+P}^M$, which have Lie products,
\[[\Phi_{O+P}^M,\Phi^{S(1+2)}]\rightarrow(\Phi_{O+P}^M\text{ or }y)\text{, with 2-forms $y\notin\Phi^{G+N(1+2)}$.}\]
For example, two members of $\Phi_O^{M(1)}$, ${\rm M}_1^{(1)}$ and ${\rm M}_{19}^{(1)}$, are
\[ {\rm M}_1^{(1)}=\frac12(-\e{19} +\e{2A} +\e{4C} -\e{7F}) \text{\ and\ }{\rm M}_{19}^{(1)} = \frac12(\e{2C} +\e{3D} +\e{4A} +\e{5B})\ec \]
which have the product $[{\rm M}_1^{(1)},{\rm M}_{19}^{(1)}] = 0$ but the second sign variation with the second and last terms negated gives the following,
\[ \big[{\rm M}_1^{(3)},{\rm M}_{19}^{(3)}\big] = \e{24} +\e{AC}\ec \]
which is also a weak isomorphism, transforming $\Phi$ to another calibration. All such products are not part of $\Phi^{S(1+2)}$.

Further, Lie products of such 2-form pairs generate all 105 2-forms of basis $\Spin{15}$ rotations, which are not part of $\Sharp{4}$. Hence the strong isomorphisms of $\Phi_{O+P}^M$ are the basis of a group distinct from $\Phi^{S(1+2)}$. The $\Cp{15}$ space, as a real module of the $\Spin{15}$ group, has been divided into the $\Phi^{S(1+2)}$ weak and strong isomorphisms on one side, and, on the other side, the mixed $\Phi_{O+P}^M$ and $\Spin{15}$ basis rotations, including derived members that transform calibrations to quasi-calibrations.
\end{proof}
\begin{definition}
  The independent basis of $\Phi_O^{C(1+2+3+4)}$ will be called $\Phi^{U(1+2+3+4)}$. The set $\Phi^{U(1)}$ contains the 14 elements of the independent basis of the 15-D $\GT$, $\Phi_O^{C(1)}$. The independent basis of the weak internal isomorphisms of $\Phi^S$ will be denoted $\Phi^T$ so that $\Phi^V=\Phi^{U(1)} +\Phi^T$ is the independent basis of $\Phi^S$. 
\end{definition}
\begin{theorem}\label{thm:6}
  The elements $\Phi_A^A+\Phi_O^{C(2+3+4)}+\Phi_P^{C(1+2+3+4)}$ are the weak internal independent elements of $\Phi^S$. The independent basis of this set, $\Phi^T$, consists of the following 91 elements, which form a basis of the 2944 calibrations of $\Harp{4}$. The algebra $\Phi^V$ is a new Lie algebra of rank 105. With $\Phi_O^{C(1)}$ provided in \tar{gt15} and denoting $\Phi_P^{C(1)}$ as $\GG{a}1,\GG{b}1,\ldots,\GG{u}1$ and $\Phi_A^{A(1)}$ as $\GA11, \GA21,\ldots,\GA71$ as well as the three sign variations of these elements, then the isomorphism basis is
  \[ \begin{split} \Phi^T = \big(&\GG{A}2, \GG{B}2, \GG{C}2, \GG{D}2, \GG{H}2,
    \GG{K}2, \GG{M}2, \GG{A}3, \GG{B}3, \GG{C}3, \GG{D}3, \GG{E}3, \GG{F}3, 
    \GG{G}3,\\
    &\GG{A}4, \GG{B}4, \GG{C}4, \GG{D}4, \GG{E}4, \GG{F}4, \GG{G}4, \GG{H}4,
    \GG{I}4, \GG{J}4, \GG{K}4, \GG{L}4, \GG{M}4, \GG{N}4, \\
    &\GG{A}4+\GG{H}4, \GG{B}4+\GG{I}4, \GG{C}4+\GG{J}4, \GG{D}4+\GG{K}4, \GG{F}4+\GG{M}4, \GG{G}4+\GG{N}4,\\
    &\GG{L}4-\GG{E}4,
    \GG{a}1, \GG{b}1, \GG{c}1, \GG{d}1, \GG{e}1, \GG{f}1, \GG{g}1, \GG{j}1, 
    \GG{m}1, \GG{o}1, \GG{p}1, \GG{q}1, \\
    &\GG{r}1, \GG{s}1, \GG{t}1, 
    \GG{u}1,
    \GG{a}2, \GG{b}2, \GG{c}2, \GG{e}2, \GG{f}2, \GG{g}2, \GG{h}2, \GG{i}2,
    \GG{k}2, \GG{l}2, \GG{m}2, \\
    &\GG{n}2,  
    \GG{a}3, \GG{b}3, \GG{c}3, \GG{d}3, \GG{e}3, \GG{f}3, \GG{g}3, \GG{j}3,
    \GG{m}3, \GG{n}3, \GG{o}3, \GG{p}3, \GG{q}3, \GG{r}3, \\
    &\GG{a}4, \GG{b}4, \GG{e}4, \GG{f}4, \GG{g}4, \GG{l}4, \GG{m}4, 
    \GA11, \GA21, \GA31, \GA41, \GA52, \GA63, \GA74\big)\es
\end{split}\]
\end{theorem}
\begin{proof}
Transforming associative calibration $\Theta$ with elements of $\Phi^U$ swaps pairs of $\Phi_O$ and pairs of $\Phi_P$ terms of $\Phi$, according to \tar{aut}. The same action occurs to the associated Fano plane algebras in \tar{fvo}. The sign variations swap the same terms but will also change pairs of signs, which changes the representation of the transformed $\Theta$. The $\Phi_A^A$ and $\Phi_P^C$ act in the same way for all sign variations. It is easy to check that $\Rot{\GG{E}2}\Theta\Rot{-\GG{E}2}=\Rot{\GG{C}2}\Theta\Rot{-\GG{C}2}$ so that the $\GG{E}2$ is dependent and the prior members of $\Phi^T$ are unique. This process is repeated on for all 168 weak internal members of $\Phi^S$ keeping the first unique calibration found to obtain the 91 independent members of $\Phi^T$. Generating $\Pin7$ reflections and pseudo-reflections multiplies these 92 unique calibrations by a cycle of degree 32, which are also all independent. Thus there are $92\times32=2944$ unique internal calibrations. The Lie product table $[\Phi^T,\Phi^T]$ includes elements of $\Phi^S$ not in $\Phi^T$ including $\Phi^{U(1)}$ so $\Phi^V$ is the independent basis of $\Phi^{S(1+2)}$. \thr{5} provides closure and members of $\Phi^V$ satisfy anticommutativity and the Jacobi identity.
\end{proof}
The action of the weak internal isomorphisms on the subalgebras of \tar{fvo} follows that of \tar{aut} apart from some Pin reflections of other subalgebras. For example, $\Phi_O^{C(1)}$ swaps $\Theta_i\leftrightarrow\Theta_j$, for $(i,j)\in\{(4,7),(5,6),(11,14),(12,13)\}$, leaving the calibration invariant, while sign variation $\Phi_O^{C(3)}$ has the same action but includes basis reflections on these elements and sign variations $\Phi_O^{C(2+4)}$ extends this to basis reflections on all subalgebras in \tar{fvo}.

Note that the dependent elements of $\Phi_P^{C(1)}$ are $\GG{h}1=\GG{b}1+\GG{L}1$, $\GG{i}1=\GG{f}1+\GG{N}1$, $\GG{k}1=\GG{e}1-\GG{K}1$, $\GG{l}1=\GG{c}1-\GG{F}1-\GG{M}1$ and $\GG{n}1=\GG{c}1-\GG{M}1$, hence these do not occur in $\Phi^T$.
\begin{theorem}\label{thm:7} The automorphisms of the set $\Harp{4}$ are
\[\Aut{\Harp{4}}=\Phi^{U(1)}\es\]
\end{theorem}
\begin{proof}
This is a graph automorphism group since it is the relative orientations of the simplex that is being preserved, not the product table derived from the calibration. The automorphism properties in \eqn{aut} are satisfied for all $\Spin{n}$ basis elements, $\phi_{\e{ij}}(y)=\Rot{\e{ij}}y\Rot{\e{ji}}$, $i,j\in\N_1^n$, $i\ne j$ and $y\in\Cl{n}$. The Spin elements form a group and are just a change of labels, which is an isomorphism but the corollary to \thr{1} shows some $\Spin{15}$ basis elements map the calibration to a quasi-calibration so can not form a group as automorphisms. \thr{5} showed that $\Psi^-\beta'=\beta'$ for $\beta\in\Phi_{A+O+P}^A$. Hence the action of all elements of $\beta\in\Phi^A$ are invariant on the 7-form calibration as full Spin rotations, $\Rot{\beta}\Phi\Rot{\beta}^{-1}=\Phi$. The proof of \ler{2} shows $\Phi^S$ is closed under the Lie product. Now, \tar{aut} shows that $\Phi^U$ swaps pairs of $\Phi_O$ and pairs of $\Phi_P$ but applying the rotations to \tar{fvo} shows that $\Phi^{C(1)}$ leaves the set of Fano planes invariant while all other internal isomorphisms of $\Phi^S$ change pairs of Fano planes. 
This means the independent basis $\Phi^U$ rotates within each plane in a way that keeps the same graph structure since rotations are pairs of reflections. This is not the case for members of $\Phi^T$.
\end{proof}
The transformations between calibrations and automorphisms of $\Harp3$ and $\Harp4$ do not have analogous structure. The $\Pin7$ reflections and pseudo-reflections provide cycles of degree 16 for any calibration of the 30 representations. For $\Pin{15}$, the reflections and pseudo-reflections provide cycles of degree 32 for the 92 unique calibrations in each primary representation. The automorphisms of $\Harp3$ and octonions are the same and due to the extra structure we expect the automorphisms of sedenions to be a subset of the automorphisms of $\Harp4$. This is analysed in the next section.

\section{Sedenion Automorphisms}
 Brown~\cite{Brown} defines three transformations that define isomorphisms of a Cayley-Dickson algebra with $N$ generators, $\AB{N}$, from those of $\AB{N-1}$. The pure transformations for $a,b\in\AB{N-1}$ and $\sigma\in\Aut{\AB{N-1}}$, are
\eqb{at15}\begin{split}
  \sigma^\prime&\maps a+b\o{N}\rightarrow a\sigma+(b\sigma)\o{N}\ec\\
  \epsilon&\maps a+b\o{N}\rightarrow a - b\o{N}\text{ and}\\
  \psi&\maps a+b\o{N}\rightarrow \frac12\Big[(-a+\beta b) +(-b+\beta a)\o{N}\Big]\es
\end{split}\eqe
The first transformation, $\sigma'$, is Schafer's automorphism, $\Phi_O^{C(1)}$. The second transformation is the conjugate of $\o{N}$, $\Phi_O^{C(2)}$, which gives, $\Phi_O^{C(1+2)}$, as
  \[\Aut{\AB{N-1}}\times\Sym2\overset{{?}}{\subseteq}\Aut{\AB{N}}\ec\]
where $\Sym2$ is the symmetric group of order 2. 
Interpreting the imaginary from Brown~\cite{Brown}, Eakin and Sathaye~\cite{Eakin} as the $\Cl{15}$ pseudoscalar, 
  \[\beta=\sqrt3\e{123456789ABCDEF}\ec\]
then $\beta\Phi_O^{(1+2)}$ forms another domain of isomorphisms. This is a result for positive definite Clifford algebra and does not apply if negative definite signature is adopted over a real field, as pointed out by Brown. Brown's transformation for pure elements can be inverse mapped into the $\Harp4$ algebra as
  \eqb{brown}\psi\maps a+b\e{u}\rightarrow\frac12(\beta-1)a -\frac12(\beta+1)b\e{u}\ec\eqe
where $u\ge8$ in hexadecimal notation up to $F$. These are the $\Cl{15}$ elements that map under \eqn{map15} to the sedenion elements that include $\o4$. The three automorphisms satisfy the properties derived by Brown,
$\psi\sigma'=\sigma'\psi$, $\epsilon\sigma'=\sigma'\epsilon$, $\epsilon^2=1$, $\psi^2=1$ and $\epsilon\psi=\psi^2\epsilon$, which defines the $\Sym3$ group, and Brown says he has proved this for $N\in\N_5^6$ also,
  \[\Aut{\AB{N-1}}\times\Sym3\overset{{?}}{\cong}\Aut{\AB{N}}\es\]
For $N=4$, this ignores the Sharp automorphisms derived from $\Phi_P$ and the other two sign variations. It also ignores the fact that the second of Brown's transformation, $\epsilon$, is a weak isomorphism. To consider this discrepancy, all the calibrations of $\Harp{4}$ are considered. This is an interpretation of the work done by Jacobson~\cite{Jacobson}.

\begin{lemma}\label{lem:3}
Internal isomorphisms of $\Phi^{S(1+2)}$, map under \eqn{map15} as \eqb{rem}\e{123456789ABCDEF}\maps\Rot{\beta}\rightarrow\pm1\text{, $\beta\in\Phi_A^A+\Phi_{O+P}^C$}.\eqe
\end{lemma}
\begin{proof}
In $\Cl7$, all $\GT$ map under $\e{1234567}$ as $\beta\rightarrow x-x=0$, $x\in\OB^0$, where $\OB^0$ is the pure basis of $\OB$. Hence $\Rot{\beta}\rightarrow 1$ in \eqn{rot7} since $\Rot{\beta}=\beta^2+\beta+1$. But the opposite parity terms, $\Phi^{G(2)}$, which are weak isomorphisms of $\Sharp3$, map as $\Rot{\beta}\rightarrow 2x$, $x\in\OB^0$, and hence are not automorphisms of $\GT$. 

For $\Spin{15}$, \eqn{rot15} maps to unity if $\e{ij}+\e{jk}+\e{mn}+\e{op}$ maps to zero or $\pm4$ and $\e{ijklmnop}$ maps to unity. This is because all terms such as $\e{ij}$ map to the same element in $\SB$, apart from sign. If the signs are all the same then the 2-form and 4-form parts of \eqn{rot15} cancel and 6-forms give $-6$. If two terms are negated then the 2-form parts cancel, the 6-form parts cancel and the 4-form parts add to $+2$. With 1 and the map of the pseudoscalar both cases give $\Rot{\beta}\rightarrow1$, for even parity 4-forms.
\end{proof}
All the terms of $\Phi^{S(1+2)}$ satisfy \eqn{rem}, as do the terms of $\Phi_O^M$. The terms of $\Phi_P^M$ would satisfy this lemma if the parity was changed but then they would not be closed under Lie products with $\Phi_O^M$. Since the mixed terms were excluded from the $\Harp4$ isomorphisms, it is found that \eqn{rem} is not reliable enough to define the automorphisms but can be used to show all $\Phi^V$ elements map to the same elements in $\SB$.

There are three candidates for the automorphisms of sedenions. 
\[\Aut{\SB}\overset{{?}}{=}\begin{cases}
  \Phi^{U(1)}\ec&\text{strong definition,} \\
  \Phi_O^{C(1+2)}+\Phi_P^{C(1+2)}\ec&\text{weak internal definition and}\\
  \Phi^V\ec&\text{weak external definition.}
\end{cases}\]
The first definition corresponds to Schafer's result and follows that of the automorphisms of $\Har{i}3$, and is analogous to $\Aut\OB$ with strong isomorphisms in both cases. \tar{op} shows the second definition forms a subalgebra of $\Phi^V$ that does not include the external isomorphisms $\Phi^{G+N}$. \ler{2} shows that including $\Phi_A$ and $\Phi_P^{C(1+2)}$ introduces $\Phi^{G+N}$ into the Lie products and \thr{5} showed that $\Rot{\beta}\Phi\Rot{-\beta}$ for $\beta\in\Phi^{G+N}$ are external calibrations. 

\begin{theorem}\label{thm:8} The sedenion automorphisms are
\[\Aut{\SB}=\Phi_{O+P}^{C(1+2)}\es\]
\end{theorem}
\begin{proof}
The proof of \thr{7} showed that members of $\Phi^V$ also satisfy the properties of automorphisms, \eqn{aut}. For sedenions there is another calibration that needs to have invariant properties under rotations. This is the 3-form associative calibration, $\Theta$ in \eqn{cal35}, with 35 terms that map the 35 quaternions in $\Cl{15}$ under \eqn{map15} to sedenions. For automorphisms of $\SB$, we expect the product table of the transformed associative calibration to match the map of the calibration itself. That is, the map $\e{123456789ABCDEF}$ of $\Rot{\beta}\Theta\Rot{-\beta}$ to equal the same map of $\Theta$, for $\beta\in\Aut{\SB}$. \ler{3} shows that this is satisfied for all internal isomorphisms of $\Phi^V$. It also shows that $\Phi^{G+N}$ do not map to identity transformations. This is an unusual situation whereby the Lie product of two candidate automorphisms can generate elements that are not automorphisms. Since this is a contradiction $\Phi^V$ is rejected, which is unfortunate because then 28 internal calibrations are not covered by the automorphisms above. These internal isomorphisms would be automorphisms if the $\Phi^{G+N}$ subalgebra of $\Phi^V$ was excluded.

The $\Pin{15}$ reflections also satisfy \eqn{aut} but do not form a group by themselves, apart from the identity and the pseudoscalar, which act as complex numbers. The action of the pseudoscalar applied to $\Phi^V$ negates the resulting calibration under conjugation, which turns all quaternions of the pure calibration into anti-quaternions. This is a weak isomorphism that negates all elements of the product table resulting in an anti-automorphism. Hence Brown's use of the pseudoscalar is rejected. The $\Spin{15}$ basis is also rejected for the same reason as in the proof of \thr{7}. This proves the theorem with independent basis $\Phi^{U(1)}+\Phi^W$ where 
\ifnum1=\AGACSE
  \[\begin{split} \Phi^W = \big(&\GG{A}2, \GG{B}2, \GG{C}2, \GG{D}2, \GG{H}2,
    \GG{K}2, \GG{M}2, \GG{a}1, \GG{b}1, \GG{c}1, \GG{d}1, \GG{e}1, \GG{f}1, \GG{g}1, \GG{j}1, \GG{m}1, \GG{o}1, \\
    &\GG{p}1, \GG{q}1, \GG{r}1, \GG{s}1, \GG{t}1, \GG{u}1,
    \GG{a}2, \GG{b}2, \GG{c}2, \GG{e}2, \GG{f}2, \GG{g}2, \GG{h}2, \GG{i}2,
    \GG{k}2, \GG{l}2, \GG{m}2, \GG{n}2\big)\es
  \end{split}\]
\else
  \[\begin{split} \Phi^W = \big(&\GG{A}2, \GG{B}2, \GG{C}2, \GG{D}2, \GG{H}2,
    \GG{K}2, \GG{M}2, \GG{a}1, \GG{b}1, \GG{c}1, \GG{d}1, \\
    &\GG{e}1, \GG{f}1, \GG{g}1, \GG{j}1, \GG{m}1, \GG{o}1, 
    \GG{p}1, \GG{q}1, \GG{r}1, \GG{s}1, \GG{t}1, \GG{u}1,\\
    &\GG{a}2, \GG{b}2, \GG{c}2, \GG{e}2, \GG{f}2, \GG{g}2, \GG{h}2, \GG{i}2,
    \GG{k}2, \GG{l}2, \GG{m}2, \GG{n}2\big)\es
  \end{split}\]
\fi
Avoiding the psuedoscalar part of Brown's solution extends his rank 21 algebra to a rank 49 algebra that covers 35 of the 91 internal calibrations of sedenions.
\end{proof}

It is easy to see how weak isomorphisms have been neglected and confused as automorphisms in $\SB$. Brown included the obvious first sign variation of $\Phi_O^C$ and $\Phi^T$ in \thr{6} shows this includes seven pure internal calibrations. Clifford algebra allows complete analysis of the isomorphisms to find 92 total internal calibrations and their associated isomorphisms that satisfy equations \eqn{aut} and form a group that include automorphisms and members that are not automorphisms. This is a unique situation and avoided here by excluding the inconsistent elements. The pseudoscalar isomorphisms introduced by Brown could be called complex automorphisms but in the Sharp algebras these are represented as pseudo-rotations and generate anti-automorphisms. The terminology weak isomorphism is thus clearer than calling them complex automorphisms in the real algebra of $\Cl{15}$ and this argument would extent to all larger Cayley-Dickson algebras. The weak parts of \thr{8} are analogous to outer automorphisms without invoking reflections. Since the sedenion case is more complicated than other automorphisms is it suggested that the definition of outer automorphisms be extended in Clifford algebra to include the weak definition. The weak isomorphisms definition is still appropriate for the internal calibrations generated by Pin reflections and the 56 members of $\Phi^T$ that are inconsistent automorphisms.

Basis rotations of $\Spin7$ were found to generate the 30 external calibrations of $\Harp3$ in~\cite{Wilmot1} when applied repeatedly. Rotating $\Phi$ or $\Theta$ repeatedly with the 105 2-forms in $\Cl{15}$ finds $64,864,800$ unique primary representations~\cite{Wilmot3}. Considering only rotations that do not generate quasi-calibrations and applying combinations of these can return the calibration back to the starting primary to find a family of 92 pure calibrations for each primary. Each of these have a 32-cycle of $\Pin{15}$ basis reflections and pseudo-reflections giving $64,864,800\times92\times32\approx191\times 10^{9}$ representations of sedenions. These calibrations classify the $\Sharp4$ ideals of $\Harp{4}$. 

The process of generating an invertible non-associative calibration from an associative calibration that maps to Cayley-Dickson algebras, $\AB{N}$, $N>3$, and generates isomorphisms and automorphisms has been demonstrated for sedenions. This process is easily extended to higher levels either by building calibrations or reverse mapping Cayley-Dickson algebras. It is expected that such calibrations will include domains that correspond to $\PB{12}$ and $\PB{14}$ power-associative rings and some rotations between the domains will generate quasi-calibrations exposing the non-associative structure of larger Cayley-Dickson algebras. There are more than 100 quasi-sedenions~\cite{Wilmot3} with three included in $\UB2$ and twelve in $\UB3$. There could be different representations in $\UB4$~\cite{Wilmot2} but the structure of Cayley-Dickson algebras should become consistent after the seven generators of $\UB4$.

The Cayley-Dickson algebra of $N$ generators is known to be a twisted mapping from matrices, $\R[\Z/2\Z]^N$~\cite{Albuquerque,Ren}. Extending this to $\Sharp{N}$ is not an untwisting as an inverse operation since this would lose the parity of the 3-cycle structure and for the sedenion case, the differentiation of the non-associative ring types shown by $\Sharp4$. For octonions that map from $8\times8$ real matrices then the map to $\Sharp3\subset\Cl7\cong\C[\Z/2\Z]^3$ is another twist. But $\Sharp{N}$ is extended to $\Harp{n}$ by multiplication of the pseudoscalar and assuming that the pseudoscalar maps to the matrix imaginary then $\Sharp3\subset\R[\Z/2\Z]^3$. Due to the Bott periodicity this extends to all $\Sharp{N}\subset\R[\Z/2\Z]^N$. All Cayley-Dickson algebras with $N>2$ generators have structure ${\rm E}_\alpha{\rm E}_\beta=-\delta_{\alpha\beta}+c_{\alpha\beta\gamma}{\rm E}_\gamma$, where $c_{\alpha\beta\gamma}$ are the antisymmetric structure constants corresponding to an associative calibration. The twisting to the Sharp algebra is given by ${\rm E}_\alpha{\rm E}_\beta=\delta_{\alpha\beta}-\abs(c_{\alpha\beta\gamma}){\rm E}_\gamma$, for $N>3$, where ${\rm E}_\alpha$ maps a single index, $\alpha$, to one of the pure $\Sharp{N}$ elements. This mapping shows that the pseudoscalar does map to the complex imaginary.

\section{Summary}
Calibrations in Clifford algebra provide a new way to analyse the structure of Cayley-Dickson algebras and this paper developed maps in order to understand this relationship. The paper also developed a new algebra series called Harp with even subalgebras called Sharp. This was defined in three ways; as a subalgebra of the Clifford algebra homomorphic to certain real matrices, as a geometric interpretation of $(N-1)$-simplices that map to $\PG{N-1}2$ and Cayley-Dickson algebras and as symmetrised Cayley-Dickson structure constants. The Sharp algebras were found to be related to ideals and the calibrations provide a classification of these ideals and hence the ideals of $\R[(\Z/2\Z)^N]$ matrices. In 7-D, the associative calibration is invertible but this is not the case in 15-D. Instead, \ler{1} was used to derive a related non-associative calibration that is invertible. Non-associative calibrations have been found in~\cite{Wilmot3} to form an infinite series and since they have the same dimensions and degrees of freedom as the associative calibration for each larger Cayley-Dickson algebra is seems likely a correspondence will always be found between the two calibrations.

The mechanism for deriving automorphisms of octonions from the calibration was extended to sedenions and found to be more complicated than that for octonions and quaternions. The Sharp algebras may be considered as duals of Cayley-Dickson algebras since they have the same ring structure of quaternions but are commuting and have positive squares, for $N>2$.
The original analysis of Schafer and Brown on sedenion automorphisms was re-analysed using Clifford algebra calibrations, which provided a clarification of these isomorphisms. The $\Sharp{N}$ algebras provide the mechanism to evaluate the automorphisms of all Cayley-Dickson algebras and this is expected to increase in complexity for the next two levels of these algebras.

The relationship of $\Sharp{N}$ to the finite geometry series $\PG{N-1}2$ was shown. Since this series is generated by projections of hyper-cubes and was discovered by Fano in 1892~\cite{Fano}, and is now found to be generated by projections of certain simplices, it is suggested that the simplex series be known as Fano hyper-volumes to distinguish it from the projective geometry of hyper-cubes. A visual representation of sedenions as the Fano volume in 3-D provides a geometric understanding of the interaction of the various non-associative rings of sedenions. This geometric approach simplified some of the proofs and it is conjectured that Fano hyper-volumes exist for all higher order Cayley-Dickson algebras.

Cawagas discovered the algebra $\PB4$~\cite{Cawagas1} and found seven copies in sedenions along with eight octonions~\cite{Cawagas2}. Using a graded notation allowed the 15 non-associative rings of sedenions to be uniquely identified and makes obvious the XOR product structure. These 15 subalgebras were inverse mapped to Fano planes in $\Cl{15}$ to find a new 7-form calibration and it is conjectured that such non-associative calibrations will exist for $\PG{N}2$, providing algebras for analysis in finite geometry. 

The algebra $\PB4$ is only one of three power-associative algebras that exist in Cayley-Dickson algebras with $\PB{12}$ and $\PB{14}$ embedded in successive doublings of sedenions. \tar{pbk} summarised how these algebras are part of a structure of six power-associative, quasi-octonion algebras and explained the naming convention. The terminology ultra-complex numbers was introduced in~\cite{Wilmot2} to distinguish power-associative algebras from complex and hyper-complex algebras and this correlates with quasi-calibrations generating quasi-algebras as well as the power-associative Cayley-Dickson algebras, all identified by the number of Types~B and~X non-associative triads.

Rotations between the 7-D $\OB$ and $\PB4$ Fano planes within the Fano volume generate quasi-calibrations that map to quasi-sedenions. Such a mapping within the sedenions would have the same effect. Hence the associative Clifford algebras have hidden knowledge of the non-associative structure of sedenions and this is expected to follow for larger Cayley-Dickson algebras in a more complicated fashion. Since Clifford algebras contain spinors and automorphisms that embed ${\rm SU}(3)$, then they may provide a sufficient description of particle physics, especially since these algebras describe Dirac algebra and could distinguish another three generations of ideals. Future work undertaken will express the Dirac and Pauli algebras as subalgebras of Sharp(3) to categorise their ideals. 

The author gratefully acknowledges the invaluable assistance and encouragement of James Chappell, Jacqui Ramagge and Derek Abbott for improving terminology, mathematical understanding and ongoing assistance with the publishing of this article. 
\ifnum0=\AGACSE
The author also acknowledges the support of an Australian Government Research Training Program Scholarship. Support from the Australian Research Council (FL240100217) is also gratefully acknowledged.
\fi
All of this work was verified with the use of a Clifford/Cayley-Dickson algebras calculator written in Python. The relevant code is the ${\it spin15\_g2.ca}$ module loaded into ${\it calcCA.py}$. The github URL for the calculator is\par\centerline{\url{https://github.com/GPWilmot/geoalg}.}

\ifnum0=\AGACSE\newpage\fi 
\section*{Appendix}
\tab{[!ht]}{auta}{Algebra $\Phi_A^{A(1)}$}
\begin{tabular}{lll}
  $\frac12(\e{89}-\e{AB}-\e{CD}+\e{EF})$ &$\frac12(\e{8A}+\e{9B}-\e{CE}-\e{DF})$ 
    &$\frac12(\e{8B}-\e{9A}-\e{CF}+\e{DE})$ \\
  $\frac12(\e{8C}+\e{9D}+\e{AE}+\e{BF})$ &$\frac12(\e{8D}-\e{9C}-\e{AF}+\e{BE})$
    &$\frac12(\e{8E}+\e{9F}-\e{AC}-\e{BD})$ \\
  $\frac12(\e{8F}-\e{9E}+\e{AD}-\e{BC})$ & &
\end{tabular}\tae
\tab{[!ht]}{aut1o}{Algebra $\Phi_O^{C(1)}$}
\setlength{\tabcolsep}{3pt}
\begin{tabular}{lll}
  $\frac12(\e{23}-\e{45}+\e{AB}-\e{CD})$ &$\frac12(-\e{13}-\e{46}-\e{9B}-\e{CE})$ 
    &$\frac12(\e{12} -\e{47} +\e{9A} -\e{CF})$\\[2pt]
  $\frac12(-\e{15}+\e{26}-\e{9D}+\e{AE})$ &$\frac12(\e{14}+\e{27}+\e{9C}+\e{AF})$ 
    &$\frac12(-\e{17}+\e{24}-\e{9F}+\e{AC})$\\[2pt]
  $\frac12(-\e{16}-\e{25}-\e{9E}-\e{AD})$ &$\frac12(\e{45}+\e{67}+\e{CD}+\e{EF})$ 
    &$\frac12(\e{46}-\e{57}+\e{CE}-\e{DF})$\\[2pt]
  $\frac12(\e{47}+\e{56}+\e{CF}+\e{DE})$ &$\frac12(-\e{26}+\e{37}-\e{AE}+\e{BF})$ 
    &$\frac12(\e{27}+\e{36}+\e{AF}+\e{BE})$\\[2pt]
  $\frac12(\e{17}+\e{35}+\e{9F}+\e{BD})$ &$\frac12(\e{25}-\e{34}+\e{AD}-\e{BC})$
    &$\frac12(\e{23}+\e{67}+\e{AB}+\e{EF})$\\[2pt]
  $\frac12(-\e{13}-\e{57}-\e{9B}-\e{DF})$ &$\frac12(\e{12}+\e{56}+\e{9A}+\e{DE})$
    &$\frac12(-\e{15}+\e{37}-\e{9D}+\e{BF})$\\[2pt]
  $\frac12(\e{24}+\e{35}+\e{AC}+\e{BD})$ &$\frac12(-\e{16}-\e{34}-\e{9E}-\e{BC})$
    &$\frac12(-\e{14}+\e{36}-\e{9C}+\e{BE})$
\end{tabular}\tae
\tab{[!ht]}{aut1p}{Algebra $\Phi_P^{C(1)}$}
\begin{tabular}{lll}
  $\frac12(\e{12}+\e{47}-\e{8B}-\e{DE})$ &$\frac12(\e{14}-\e{27}+\e{8D}-\e{BE})$
    &$\frac12(\e{17}+\e{24}+\e{8E}+\e{BD})$\\[2pt]
  $\frac12(\e{12}-\e{56}-\e{8B}+\e{CF})$ &$\frac12(-\e{15}-\e{26}+\e{8C}+\e{BF})$
    &$\frac12(\e{16}-\e{25}-\e{8F}+\e{BC})$\\[2pt]
  $\frac12(\e{13}-\e{46}+\e{8A}-\e{DF})$ &$\frac12(\e{14}+\e{36}+\e{8D}+\e{AF})$
    &$\frac12(\e{16}-\e{34}-\e{8F}+\e{AD})$\\[2pt]
  $\frac12(\e{13}-\e{57}+\e{8A}-\e{CE})$ &$\frac12(-\e{15}-\e{37}+\e{8C}+\e{AE})$
    &$\frac12(\e{17}-\e{35}+\e{8E}-\e{AC})$\\[2pt]
  $\frac12(\e{23}+\e{45}-\e{89}-\e{EF})$ &$\frac12(\e{24}-\e{35}+\e{8E}-\e{9F})$
    &$\frac12(\e{25}+\e{34}-\e{8F}-\e{9E})$\\[2pt]
  $\frac12(\e{23}-\e{67}-\e{89}+\e{CD})$ &$\frac12(-\e{26}-\e{37}-\e{8C}-\e{9D})$
    &$\frac12(\e{27}-\e{36}+\e{8D}-\e{9C})$\\[2pt]
  $\frac12(\e{45}-\e{67}-\e{89}+\e{AB})$ &$\frac12(\e{46}+\e{57}-\e{8A}-\e{9B})$
    &$\frac12(\e{47}-\e{56}-\e{8B}+\e{9A})$
\end{tabular}\tae
\tab{[!ht]}{aut2o}{Algebra $\Phi_O^{M(1)}$}
\setlength{\tabcolsep}{3pt}
\begin{tabular}{lll}
  $\frac12(-\e{19}+\e{2A}+\e{4C}-\e{7F})$ &$\frac12(-\e{1A}+\e{29}+\e{4F}-\e{7C})$
    &$\frac12(\e{1C}+\e{2F}-\e{49}-\e{7A})$\\[2pt]
  $\frac12(-\e{1F}+\e{2C}+\e{4A}-\e{79})$ &$\frac12(-\e{19}+\e{2A}-\e{5D}+\e{6E})$
    &$\frac12(\e{1A}-\e{29}+\e{5E}-\e{6D})$\\[2pt]
  $\frac12(\e{1D}+\e{2E}+\e{59}+\e{6A})$ &$\frac12(\e{1E}+\e{2D}-\e{5A}-\e{69})$
    &$\frac12(\e{19}+\e{3B}-\e{4C}-\e{6E})$\\[2pt]
  $\frac12(\e{1B}+\e{39}+\e{4E}+\e{6C})$ &$\frac12(\e{1C}-\e{3E}-\e{49}+\e{6B})$
    &$\frac12(\e{1E}+\e{3C}-\e{4B}-\e{69})$\\[2pt]
  $\frac12(\e{19}+\e{3B}+\e{5D}+\e{7F})$ &$\frac12(\e{1B}+\e{39}+\e{5F}+\e{7D})$
    &$\frac12(\e{1D}-\e{3F}+\e{59}-\e{7B})$\\[2pt]
  $\frac12(\e{1F}+\e{3D}+\e{5B}+\e{79})$ &$\frac12(\e{2A}-\e{3B}+\e{4C}-\e{5D})$
    &$\frac12(\e{2B}-\e{3A}-\e{4D}+\e{5C})$\\[2pt]
  $\frac12(\e{2C}+\e{3D}+\e{4A}+\e{5B})$ &$\frac12(\e{2D}-\e{3C}+\e{4B}-\e{5A})$
    &$\frac12(\e{2A}-\e{3B}+\e{6E}-\e{7F})$\\[2pt]
  $\frac12(\e{2B}-\e{3A}+\e{6F}-\e{7E})$ &$\frac12(\e{2E}-\e{3F}+\e{6A}-\e{7B})$
    &$\frac12(\e{2F}+\e{3E}-\e{6B}-\e{7A})$\\[2pt]
  $\frac12(\e{4C}-\e{5D}+\e{6E}-\e{7F})$ &$\frac12(\e{4D}-\e{5C}+\e{6F}-\e{7E})$
    &$\frac12(\e{4E}-\e{5F}+\e{6C}-\e{7D})$\\[2pt]
  $\frac12(\e{4F}+\e{5E}-\e{6D}-\e{7C})$ & &
\end{tabular}\tae
\tab{[!ht]}{aut2p}{Algebra $\Phi_P^{M(1)}$}
\setlength{\tabcolsep}{4pt}
\begin{tabular}{lll}
  $\frac12(\e{18}+\e{2B}+\e{4D}+\e{7E})$ &$\frac12(\e{1B}+\e{28}+\e{4E}+\e{7D})$
    &$\frac12(\e{1D}+\e{2E}+\e{48}+\e{7B})$\\[2pt]
  $\frac12(\e{1E}+\e{2D}+\e{4B}+\e{78})$ &$\frac12(\e{18}+\e{2B}+\e{5C}+\e{6F})$
    &$\frac12(\e{1B}+\e{28}+\e{5F}+\e{6C})$\\[2pt]
  $\frac12(\e{1C}+\e{2F}+\e{58}+\e{6B})$ &$\frac12(\e{1F}+\e{2C}+\e{5B}+\e{68})$
    &$\frac12(\e{18}+\e{3A}+\e{4D}+\e{6F})$\\[2pt]
  $\frac12(\e{1A}+\e{38}+\e{4F}+\e{6D})$ &$\frac12(\e{1D}+\e{3F}+\e{48}+\e{6A})$
    &$\frac12(\e{1F}+\e{3D}+\e{4A}+\e{68})$\\[2pt]
  $\frac12(\e{18}+\e{3A}+\e{5C}+\e{7E})$ &$\frac12(\e{1A}+\e{38}+\e{5E}+\e{7C})$
    &$\frac12(\e{1C}+\e{3E}+\e{58}+\e{7A})$\\[2pt]
  $\frac12(\e{1E}+\e{3C}+\e{5A}+\e{78})$ &$\frac12(\e{28}+\e{39}+\e{4E}+\e{5F})$
    &$\frac12(\e{29}+\e{38}+\e{4F}+\e{5E})$\\[2pt]
  $\frac12(\e{2E}+\e{3F}+\e{48}+\e{59})$ &$\frac12(\e{2F}+\e{3E}+\e{49}+\e{58})$
    &$\frac12(\e{28}+\e{39}+\e{6C}+\e{7D})$\\[2pt]
  $\frac12(\e{29}+\e{38}+\e{6D}+\e{7C})$ &$\frac12(\e{2C}+\e{3D}+\e{68}+\e{79})$
    &$\frac12(\e{2D}+\e{3C}+\e{69}+\e{78})$\\[2pt]
  $\frac12(\e{48}+\e{59}+\e{6A}+\e{7B})$ &$\frac12(\e{49}+\e{58}+\e{6B}+\e{7A})$
    &$\frac12(\e{4A}+\e{5B}+\e{68}+\e{79})$\\[2pt]
  $\frac12(\e{4B}+\e{5A}+\e{69}+\e{78})$ &&
\end{tabular}\tae

\tab{[ht]}{autg}{Algebra $\Phi^{G(1)}$}
\begin{tabular}{ccccc}
 $\frac12(\e{23} -\e{45})$ &$\frac12(-\e{13}-\e{46})$ &$\frac12(\e{12}-\e{47})$ &$\frac12(-\e{15}+\e{26})$ \\[2pt]
 $\frac12(\e{14}+\e{27})$ &$\frac12(-\e{17}+\e{24})$ &$\frac12(-\e{16}-\e{25})$ &$\frac12(\e{45}+\e{67})$ \\[2pt]
 $\frac12(\e{46}-\e{57})$ &$\frac12(\e{47}+\e{56})$ &$\frac12(-\e{26}+\e{37})$ &$\frac12(\e{27}+\e{36})$ \\[2pt]
 $\frac12(\e{17}+\e{35})$ &$\frac12(\e{25}-\e{34})$ &$\frac12(\e{23} +\e{67})$ &$\frac12(-\e{13} -\e{57})$ \\[2pt]
 $\frac12(\e{12} +\e{56})$ &$\frac12(-\e{15} +\e{37})$ &$\frac12(\e{24} +\e{35})$ &$\frac12(-\e{16} -\e{34})$ \\[2pt]
 $\frac12(-\e{14} +\e{36})$
\end{tabular}\tae

\tab{[ht]}{autn}{Algebra $\Phi^{N(1)}$}
\begin{tabular}{ccccc}
  $\frac12(\e{8B}+\e{DE})$ &$\frac12(-\e{9E}+\e{AD})$ &$\frac12(-\e{8A}+\e{9B})$ &$\frac12(\e{8F}+\e{AD})$ \\[2pt]
  $\frac12(\e{9E}-\e{BC})$ &$\frac12(\e{AB}-\e{EF})$ &$\frac12(-\e{CD}+\e{EF})$  &$\frac12(\e{8E}-\e{AC})$ \\[2pt]
  $\frac12(-\e{9D}+\e{BF})$ &$\frac12(\e{AE}+\e{BF})$ &$\frac12(\e{89}+\e{CD})$ &$\frac12(\e{8D}+\e{9C})$ \\[2pt]
  $\frac12(-\e{8C}+\e{BF})$ &$\frac12(\e{8B}+\e{CF})$ &$\frac12(-\e{8E}+\e{AC})$ &$\frac12(\e{9C}+\e{AF})$ \\[2pt]
  $\frac12(\e{9F}+\e{BD})$ &$\frac12(\e{8A}+\e{DF})$ &$\frac12(\e{8C}+\e{AE})$ &$\frac12(\e{9A}-\e{CF})$ \\[2pt]
  $\frac12(-\e{8E}+\e{9F})$ &$\frac12(\e{8A}-\e{CE})$ &$\frac12(\e{AD}-\e{BC})$ &$\frac12(\e{9B}+\e{DF})$ \\[2pt]
  $\frac12(\e{8D}-\e{AF})$ &$\frac12(\e{9C}+\e{BE})$ &$\frac12(\e{AC}+\e{BD})$ &$\frac12(\e{9A}-\e{DE})$ \\[2pt]
  $\frac12(\e{CE}+\e{DF})$ &$\frac12(\e{8F}+\e{BC})$ &$\frac12(\e{89}+\e{AB})$ &$\frac12(\e{AB}-\e{CD})$ \\[2pt]
  $\frac12(\e{9D}+\e{AE})$ &$\frac12(\e{8F}+\e{9E})$ &$\frac12(\e{89}+\e{EF})$ &$\frac12(-\e{9B}+\e{CE})$ \\[2pt]
  $\frac12(\e{BC}-\e{9D})$ &$\frac12(-\e{8D}+\e{AF})$ &$\frac12(-\e{CF}+\e{DE})$ &$\frac12(\e{8E}+\e{BD})$ \\[2pt]
  $\frac12(\e{8D}-\e{BE})$ &$\frac12(\e{9F}-\e{AC})$ &$\frac12(\e{AF}-\e{BE})$ &$\frac12(\e{8B}+\e{9A})$ \\
\end{tabular}
\ifnum0=\AGACSE\vskip 20ex\fi
\tae

\begin{thebibliography}{24}
\bibitem{Morano} G.~Moreno, The zero divisors of the Cayley–Dickson algebras over the real numbers, Bol. Soc. Mat. Mex. (tercera serie) {\bf 4:1}, 1998, pp.~13-28.
\bibitem{Zhilina} S.~A.~Zhilina, Orthogonality graphs of real Cayley-Dickson algebras. Part I: Doubly alternative zero divisors and their hexagons. Int.~J.~Algebra Comput. {\bf 31(4)}, 2021, pp. 663--689.
\bibitem{Wilmot2} G.~P.~Wilmot, Structure of the Cayley-Dickson algebras, 2025, arXiv:2505.11747 [math.RA]. 
\bibitem{Dray} T.~Dray and C.~A.~Manogue, The geometry of the octonions. World Scientific Publishing Co. Pte. Ltd, Singapore, 2015, ISBN 978-9814401814.
\bibitem{Agricola} I.~Agricola, Old and new on the exceptional
group $\GT$, Notices of the AMS, {\bf 55}, No.~8, (2008) pp.~922--929, MR 2441524
\bibitem{Schafer} R.~D.~Schafer, On the Algebras Formed by the Cayley-Dickson Process, American Journal of Mathematics, {\bf 76:2}, 1954, pp.~435-446.
\bibitem{Brown} R.~B.~Brown, On generalized Cayley-Dickson algebras, Pacific Journal of Mathematics, {\bf 20}, 1967, pp.~415--422.
\bibitem{Harvey} F.~Reese Harvey, {\it Spinors and Calibrations}, Perspectives in Mathematics, Vol.~9. Academic Press, Inc. 1990.
\bibitem{Porteous} I.~R.~Porteous, {\it Topological Geometry}, Cambridge U. Press, 1981.
\bibitem{Bryant} R.~L.~Bryant, Metrics with exceptional holonomy. Ann. of Math. {\bf 126}, 1987, pp.~525--576.
\bibitem{Wilmot1} G.~P.~Wilmot, Construction of Exceptional Lie Algebra G2 and Non-associative Algebras Using Clifford Algebra, Adv. Appl. Clifford Algebras, {\bf 36:24}, 2026, https://doi.org/10.1007/s00006-025-01423-5.
\bibitem{Cawagas1} R.~E.~Cawagas, On the structure and zero divisors of the Cayley-Dickson sedenion algebra, Discussiones Mathematicae, General Algebra and Applications 24, 2004, pp.~251--265.
\bibitem{Fano} Fano, G. Sui postulati fondamentali della geometria proiettiva in uno spazio lineare a un numero qualunque di dimensioni (On the fundamental postulates of perspective geometry in a linear space of any number of dimensions), Giornale di Matematiche, {\bf 30}, 1892, pp.~106--132.
\bibitem{Cawagas2} R.~E.~Cawagas, The subloop structure of the Cayley-Dickson sedenion loop, Matimyas Matematika, {\bf 28} (2005), pp.~10--20.
\bibitem{Polster} B.~Polster, A Geometric Picture Book, Springer Science+Business Media, New York (1998) ISBN~978-1-4612-6426-2
\bibitem{Guy} R.~K.~Guy, The Unity of Combinatorics, Combinatorics Advances, Ed.~C.~J.~Colbourn, E.~S.~Mahmoodian, pp.~129--159,
Kluwer Academic Publishers, 1995.
\bibitem{Saniga} M.~Saniga, F.~Holweck, P.~Pracna, From Cayley-Dickson Algebras to Combinatorial Grassmannians, Mathematics 2015, 3, pp.~1192--1221, doi:10.3390/math3041192.
\bibitem{Gresnigt} N.~G.~Gresnigt, A sedenion algebraic representation of three colored fermion generations, J. Phys.: Conf. Ser. {\bf 2667}, 012061, 2023, doi:10.1088/1742-6596/2667/1/012061.
\bibitem{Lounesto} P.~Lounesto, Clifford Algebras and Spinors, Clifford Algebras and Their Applications in Mathematical Physics, edited by J.~S.~R.~Chisholm and A.~K.~Common, D.~Reidel Publishing Company, 1986, pp.~25--37.
\bibitem{Wilmot3} G.~P.~Wilmot, Sedenion and Ultracomplex Symmetric Cryptography, 2026, arXiv:2608.00636 [math.RA].
\bibitem{Eakin} P.~Eakin and A.~Sathaye, On automorphisms and derivations of Cayley-Dickson algebras, J.~Algebra, {\bf 129(2)}, 1990, pp.~263--278.
\bibitem{Jacobson} N.~Jacobson, Composition algebras and their automorphisms, Rend. Circ. Mat. Palermo {\bf 2:7}, 1958, pp.~55--80.
\bibitem{Albuquerque} H.~Albuquerque and S.~Majid, Quasialgebra structure of the octonions, Journal of Algebra, {\bf 220:1}, 1999, pp.~188--224, doi:10.1006/jabr.1998.7850.
\bibitem{Ren} G.~Ren and X.~Zhao, The Explicit Twisted Group Algebra Structure of the Cayley–Dickson Algebra, Advances in Applied Clifford Algebras, {\bf 33:49}, 2023, doi:10.1007/s00006-023-01296-6.
\end{thebibliography}
\end{document}